\documentclass[11pt,a4paper,amssymb,amsmath, tightenlines]{article}

\usepackage{amsmath, amssymb, amsthm, latexsym, mathrsfs, url, xcolor, epsfig, graphics, float, setspace}
\usepackage{sectsty}
\usepackage[T1]{fontenc}
\numberwithin{equation}{section}

\newtheorem{prop}{Proposition}[section]
\newtheorem{lem}{Lemma}[section]
\newtheorem{cor}{Corollary}[section]
\theoremstyle{definition}

\theoremstyle{remark}
\newtheorem{rem}{Remark}[section]
\newtheorem{ex}{Example}[section]

\newcommand{\e}{{\rm e}}
\newcommand{\E}{\mathbb{E}}

\newcommand{\F}{\mathcal{F}}

\newcommand{\PR}{\mathbb{P}}

\newcommand{\indi}[1]{1\hspace{-.09cm}\textup{\textrm{l}}}

\newcommand{\nn}{\nonumber}

\begin{document}
\title{\bf Stochastic Modelling with \\ Randomised Markov Bridges}
\author{Andrea Macrina{$^{\star\,\ddag}$}\footnote{Corresponding author: a.macrina@ucl.ac.uk}\,\, and Jun Sekine{\,$^{\S}$} 
\\ \\ 
{$^{\star}$Department of Mathematics, University College London} \\ 
{London WC1E 6BT, United Kingdom} \\ 
{$^{\ddag}$African Institute for Financial Markets and Risk Management} \\ {University of Cape Town} 
\\ {Rondebosch 7701, South Africa} \\ 
{$^{\S}$Division of Mathematical Science for Social Systems} \\ 
{Graduate School of Engineering Science, Osaka University} \\ 
{ 1-3 Machikaneyama, Toyonaka, Osaka 560-8531, Japan}}
\date{13 December 2019}
\maketitle
\vspace{-0.75cm}
\begin{abstract}
\noindent
We consider the filtering problem of 
estimating a hidden random variable $X$ by noisy observations.
The noisy observation process is constructed 
by a randomised Markov bridge (RMB) $(Z_t)_{t\in [0,T]}$ 
of which terminal value is set to $Z_T=X$. 
That is, at the terminal time $T$, the noise of the bridge process 
vanishes and the hidden random variable $X$ is revealed. 
We derive the explicit filtering formula, governing the dynamics of the conditional probability process, 
for a general RMB. It turns out that the conditional probability 
is given by a function of current time $t$, the current observation $Z_t$, 
the initial observation $Z_0$, 
and the {\it a priori} distribution $\nu$ of $X$ at $t=0$. As an example for an RMB we explicitly construct the skew-normal randomised diffusion bridge and show how it can be utilised to extend well-known commodity pricing models and how one may propose novel stochastic price models for financial instruments linked to greenhouse gas emissions.  
\\\vspace{-0.2cm}
\\\vspace{-0.2cm}\\
\end{abstract}

\section{Introduction}
Let $E \subset {\mathbb R}^n$ 
and let $T\in (0,\infty)$. 
On a probability space 
$(\Omega,{\mathcal F}, {\mathbb P})$, 
consider an $E$-valued 
$(y,T,z)$-Markov bridge $Y^{(y,T,z)}:=({Y}^{(y,T,z)}_t)_{t\in [0,T]}$
so that
\begin{equation}
 Y^{(y,T,z)}_0=y,\quad Y^{(y,T,z)}_T=z 
\quad\text{a.s.,}
\end{equation}
and an $E$-valued random variable $X$, 
which is independent of $Y^{(y,T,z)}$. 
Here, 
by a $(y,T,z)$-Markov bridge, we mean 
a process obtained by conditioning a Markov process
$Y:=(Y_t)_{t\ge 0}$ 
to start in $y\in E$ at time $0$ and arrive at 
$z\in E$ at time $T\in (0,\infty)$.
For its construction, we follow Fitzsimmons et al. (1993). 
We define the process $Z:=(Z_t)_{t\in [0,T]}$ by
\[
 Z:=Y^{(y,T,X)},
\]
which we call {\it randomised Markov bridge} (RMB).
We further define 
\[
 {\mathcal F}^Z_t:=\sigma(Z_{s}; s\in [0,t]),
\quad t\in [0,T], 
\]
and let $Z$ be the noisy observation process
of the hidden random variable $X$.
At the terminal time $T$, 
the noise in the bridge process $Y^{(y,T,z)}$ vanishes
and the hidden variable $X(=Z_T)$ is revealed. 
We are interested in the stochastic filtering 
problem of estimating the hidden random variable $X$ 
through the observation of $Z$. That is, 
we are interested in computing the conditional probability,
\begin{equation}
 \pi_t(dx):={\mathbb P}\left( X\in dx | {\mathcal F}^Z_t\right),
\quad t\in [0,T),
\end{equation}
and the conditional expectation,
\begin{equation}
 \pi_t(f):= {\mathbb E}\left[ f(X) | {\mathcal F}^Z_t\right]
=\int_E f(x)\pi_t(dx),
\quad t\in [0,T),
\end{equation}
where $f:E\to {\mathbb R}$ is Borel-measurable 
so that $f(X)$ is integrable.

The problem (1.2)-(1.3) of filtering a hidden random variable $X$ 
by observing a sub-class of RMB processes $Z$ 
has been studied in a financial context by, e.g., 
Brody {\it et al.} (2008), 
Hoyle {\it et al.} (2011), and 
Filipovi\'c {\it et al.} (2012). 
In those financial applications, $f(X)$ represents the cash flow of 
a financial asset that is paid at the terminal date $T$, and 
$Z$ is called the information process. 
The conditional expectation (1.3) is used to model the asset price process
$(S_t)_{t\in [0,T)}$ by
\begin{equation}\label{pricing_formula}
S_t:={\mathrm e}^{-r(T-t)} \pi_t(f)
={\mathbb E}\left[ {\mathrm e}^{-r(T-t)}f(X)\bigm|{\mathcal F}^Z_t\right],
\end{equation}
where $r\ge 0$ is the constant risk-free interest rate
and ${\mathbb P}$ is regarded as being the so-called 
risk-neutral probability measure.
In the above-mentioned works, 
Brownian random bridges and, more generally, 
L\'evy random bridges (LRB) are employed,  
and the stochastic dynamics of asset prices are derived. 

Randomised Markov bridges generalise the class of information processes 
which can be applied to develop information-based asset pricing models. 
Such models describe the stochastic nature of asset prices as 
the market's price adjustments to the information flow about market
factors, 
modelled by the hidden random variable $X$, 
which affect the payoff function of a financial contract. 
In this setting, the effect of emerging information about market factors 
is manifested in the conditional probability process underlying 
the stochastic dynamics of asset prices. 
By extending the class of information processes with RMBs, 
the class of information-based asset price models is also extended. 
The generalised class can now be constructed by starting off 
with any Markov process, also those that have dependent increments, 
leading to information processes beyond the LRB class 
(which are constructed on the basis of L\'evy processes). 
Further work making use of LRBs, that now may be extended by RMBs, 
includes Macrina (2014) on heat kernel asset pricing models, 
Hoyle {\it et al.} (2015) on insurance claim reserving, 
and Cr\'epey {\it et al.} (2016) on interest rate modelling.

The main contribution of the present article is 
to provide the explicit representations of (1.2) and (1.3) 
for a {\it general} RMB,
focusing on its interesting features from a stochastic filtering viewpoint.
It is well-known that
in the general stochastic filtering problem, 
the conditional probability (1.2) 
has an infinite-dimensional structure; 
it is the solution to the measure-valued Kushner-Stratonovich 
stochastic differential equation (SDE). 
In the problem we consider, unlike in the general situation, 
we obtain the relation
\begin{equation}
{\mathbb P}\left( Z_t \in dy | {\mathcal F}^Z_s\right) 
={\mathbb P}\left( Z_t \in dy | Z_s, Z_0 \right)
=: P_{s,t}(Z_s, dy| Z_0)
\end{equation}
for $0\le s< t\le T$ (see Propositions 2.1 and 2.2, 
and Remark 2.3 for details). 
In words, 
we can regard $Z$ as satisfying the Markov property 
with respect to its natural filtration 
$({\mathcal F}^Z_t)_{t\in [0,T)}$, 
once the initial value $Z_0$ is fixed (see Proposition \ref{RMB-transdens}).
From the above relation, it follows that
\[
\pi_t(dy)={\mathbb P}\left(Z_T\in dy| {\mathcal F}^Z_t\right)
=P_{t,T}(Z_t, dy| Z_0).
\]
We thus observe that
the pair of observations $(Z_0,Z_t)$ determines $\pi_t(dz)$, 
and the past observation $(Z_s)_{s\in (0,t)}$
is not necessary for the computation of the conditional probability. 
As a consequence, the dynamics of 
$(\pi_t(dz))_{t\in [0,T)}$ can be determined
by a finite-dimensional Markovian SDE, see Proposition \ref{prop3}.
\begin{rem}
We refer to Baudoin (2002) 
and \c{C}etin and Danilova (2016)
for related pieces of work: conditioned stochastic differential equations (CSDE) introduced by Baudoin (2002), and 
the weak conditioning of (Markovian) SDEs via $h$-transforms
considered by \c{C}etin \& Danilova (2016) 
are directly related with our RMBs 
(with respect to the filtration $({\mathcal F}^Z_t)_{t\in [0,T)}$).
Although the works by Baudoin (2002), 
\c{C}etin and Danilova (2016), 
and the analysis presented in this article overlap in places, 
the following viewpoint and motivation appear to be different. 
For example, (a) we are interested in providing 
a stochastic filtering interpretation of RMBs, 
see Propositions \ref{prop1} and \ref{prop3},
and (b) we consider general RMBs whereas 
Baudoin (2002) and \c{C}etin \& Danilova (2016)
study SDEs driven by Brownian motions only.
\end{rem}
In the following sections, 
after preparing the setup in detail, 
we state our results, give additional explanations in the remarks, and
provide the proofs. In Section 4, we introduce a new RMB, the {\it skew-normal randomised diffusion bridge} (SNRDB) which is applied to commodity pricing, and which, in doing so, extends the models by Gibson \& Schwarz (1990) and Schwarz (1997). In our view we also find a natural way to apply this class of pricing models to securitise climate risk associated to greenhouse gas emissions. 

\section{Setup and results}
Let $E \subset {\mathbb R}^n$ be a Borel state space, 
and let $T\in (0,\infty)$ be a given constant. 
For the construction of the $E$-valued 
$(y,T,z)$-Markov bridge $Y^{(y,T,z)}:=({Y}^{(y,T,z)}_t)_{t\in [0,T]}$
which satisfies (1.1), 
we follow Section 2 of Fitzsimmons et al. (1993). 
We consider an $E$-valued strong Markov process 
$Y:=(Y_t)_{t\in [0,T]}$ with c\`adl\`ag sample paths, 
which is realised as the coordinate process $Y$
on $\Omega^1_T$, that is, on the space of right-continuous paths
from $[0,T)$ to $E$  that have left limits on $(0,T)$. 
The law of the Markov process $Y$ that starts from $y$ 
is denoted by $\tilde{\mathbb P}^1_y$, 
and the natural filtration of $Y$ is denoted by
$({\mathcal F}^1_t)_{t\in [0,T]}$.
We assume that the transition probability of the Markov process $Y$ 
has the density
\[
 \tilde{P}_t(x,dy)
=\tilde{p}_t(x,y)m(dy)
\]
with respect to a $\sigma$-finite measure $m$ on $E$
such that the Chapman-Kolmogorov identity
\[
 \tilde{p}_{t+s}(x,z)=\int_{E} \tilde{p}_t(x,y) \tilde{p}_s(y,z)m(dy)
\]
holds true, where $t>0$ and $s>0$ ($s+t\le T$), $x\in E$ and $z\in E$, and 
$\tilde{p}_t(x,y)>0$ for all $(t,x,y)\in (0,T]\times E\times E$. 
Moreover, we assume the following duality: 
Associated with $(Y,(\tilde{P}_t)_{t\ge 0})$,
there exists $\bigl( \hat{Y},(\hat{P}_t)_{t\ge 0} \bigr)$, 
where $\hat{Y}$ is the second Markov process with the state space $E$
and $(\hat{P}_t)_{t\ge 0}$ is the transition probabilities of $\hat{Y}$, which 
satisfies
\[
 \int_E (\tilde{P}_t f_1)(x) f_2(x)m(dx)
=\int_E f_1(x) (\hat{P}_tf_2)(x)m(dx)
\]
for $t\ge 0$ and bounded Borel measurable $f_1,f_2: E \to {\mathbb R}$, 
where we write
$(\tilde{P}_t f_1)(x):= \int_E f_1(y)\tilde{P}_t(x,dy)$
and 
$(\hat{P}_t f_2)(x):= \int_E f_2(y)\hat{P}_t(x,dy)$.
Then, by Proposition 1 of Fitzsimmons et al. (1993), 
we see the following:
\begin{itemize}
 \item[(a)] Let 
${\mathcal F}^1_{T-}:=
\sigma\left( \cup_{0\le t<T} {\mathcal F}^1_t \right)$. 
For each $y$ and $z\in E$, 
we can construct the probability measure ${\mathbb P}^1_{y,z}$
on $(\Omega^1_T, {\mathcal F}^1_{T-})$  that satisfies
\[
d{\mathbb P}^1_{y,z} \bigm|_{{\mathcal F}^1_t}
=
\Lambda_t^{(y,T,z)}
d\tilde{\mathbb P}^1_{y} \bigm|_{{\mathcal F}^1_t}
\]
for all $t\in [0,T)$, where
\[
\Lambda_t^{(y,T,z)}:=
\frac{\tilde{p}_{T-t}(Y_t,z)}{\tilde{p}_T(y,z)}.
\]
 \item[(b)] The law of $(Y_t)_{t\in [0,T)}$ under ${\mathbb P}^1_{y,z}$
is equal to that of the $(y,T,z)$-Markov bridge. In particular, 
it holds that
\[
 {\mathbb P}^1_{y,z}\left( Y_0=y, Y_{T-}=z\right)=1.
\]
The corresponding transition densities that satisfy
\[
 p^{(z,T)}(s,y;t,y')m(dy')
:={\mathbb P}^{1}_{y,z}
\left( Y_t\in dy'| Y_s=y\right)
\]
for $y,y',z\in E$ and $0< s< t< T$
are expressed by
\[
 p^{(z,T)}(s,y;t,y')
=\frac{\tilde{p}_{t-s}(y,y') \tilde{p}_{T-t}(y',z)}
{\tilde{p}_{T-s}(y,z)}.
\]
\end{itemize}
Furthermore, we introduce another probability space 
$(\Omega_2,{\mathcal F}_2,{\mathbb P}_2)$ on which we consider the random variable $X$ with law
$\nu:={\mathbb P}_2\circ X^{-1}$.
Let $\Omega:=\Omega^1_T \times \Omega_2$,
${\mathcal F}_t:= {\mathcal F}^1_t \otimes {\mathcal F}_2$, 
${\mathcal F}_{T-}={\mathcal F}^1_{T-} \otimes {\mathcal F}_2$, 
and let ${\mathbb P}_y:={\mathbb P}^1_y \otimes {\mathbb P}_2$ such that
\begin{equation}
d{\mathbb P}_{y}
 (\omega_1,\omega_2)\bigm|_{{\mathcal F}_t} 
:= {\Lambda}_t^{(y,T,X)}(\omega_1,\omega_2)
d\tilde{\mathbb P}^1_y(\omega_1) \otimes d{\mathbb P}_2(\omega_2) 
\bigm|_{{\mathcal F}_t}
\end{equation}
is satisfied for $t\in [0,T)$ where 
\begin{equation}
\Lambda_t^{(y,T,X)}(\omega_1,\omega_2)
:=\frac{\tilde{p}_{T-t}\left( Y_t(\omega_1), X(\omega_2)\right)}
{\tilde{p}_T(y,X(\omega_2))}.
\end{equation}
Then, on the filtered product space 
$(\Omega,{\mathcal F}_{T-}, {\mathbb P}_{y}, 
({\mathcal F}_t)_{t\in [0,T)})$, 
the RMB $(Z_t)_{t\in [0,T]}$ 
is obtained by setting 
$Z_t(\omega_1,\omega_2):=Y_t(\omega_1)$ for $t\in [0,T)$ 
and $Z_T(\omega_1,\omega_2):=X(\omega_2)$.
\begin{rem}
Intuitively, the RMB is constructed as $Z:=Y^{(y,T,X)}$, 
by inserting the random variable $X$
into the terminal value of the Markov bridge $Y^{(y,T,\cdot)}$. 
For the validity of such a direct approach,  
the measurability of the function $z\mapsto Y^{(y,T,z)}$ is necessary.
The above construction of an RMB avoids this measurability issue. 
However, there are situations where the measurability issue 
can be resolved directly:
For example, if we consider the SDE for $Y^{(y,T,z)}$, 
starting from the Markovian SDE (2.4) below, 
then, as we will see in (2.6), we have
\begin{multline*}
 dY^{(y,T,z)}_t = 
\left\{ b\left( Y^{(y,T,z)}_t\right)+ 
(\sigma\sigma^\top)\left( Y^{(y,T,z)}_t\right)
\nabla \ln \tilde{p}_{T-t}\left( Y^{(y,T,z)}_t,z\right)
\right\} dt \\
+\sigma\left( Y^{(y,T,z)}_t\right)dW_t
\quad 
\left(Y^{(y,T,z)}_0=y\right).
\end{multline*}
The continuity of $z\mapsto Y^{(y,T,z)}$ 
can be shown by discussing the continuity of the solution to the above SDE
with respect to the parameter $z$, 
which is a standard result in SDE theory
if we assume the continuity of $z\mapsto \tilde{p}_t(y,z)$.
\end{rem}
\begin{rem}
For the construction of Markovian bridges, 
we follow classical results obtained in Fitzsimmons et al. (1993). 
There are several recent studies, where extended results are discussed:
For example, we refer to Chaumont \& Uribe Bravo (2011), 
and \c{C}etin \& Danilova (2016). 
Both studies succeed in constructing Markov bridges
without assuming the duality condition of the underlying Markov
 process. The former work derives and utilises 
the weak continuity property of a Markov
 bridge with respect to its terminal value. 
The latter work is interested in SDE representations of the bridges
driven by Brownian motions. 
Interestingly, the weak conditioning discussed in the latter work 
directly links to an RMB with respect to its natural filtration
$({\mathcal F}^Z_t)_{t\in [0,T)}$; see Remark 2.6 below.
\end{rem}

\subsection{Filtering results}
For the filtering problem, we have the following:
\begin{prop}\label{prop1}
For $t\in [0,T)$, the conditional probability (1.2)
is given by
\[
\pi_t(dz)
= 
\frac{\displaystyle \frac{\tilde{p}_{T-t}(Z_t, z)}{\tilde{p}_T(Z_0,z)} \nu (dz)}
{\displaystyle \int_E \frac{\tilde{p}_{T-t}(Z_t, z')}
{\tilde{p}_T(Z_0,z')} \nu (dz')}.
\]
\end{prop}

\begin{rem}
From a stochastic filtering viewpoint, 
it is interesting that
the conditional probability $\pi_t(dz)$
is expressed by a function of $t$, $Z_t$, $Z_0$, and $\nu$. 
At each point in time, 
$\pi_t(dz)$ can be computed by inserting the current observation $Z_t$ 
and the initial observation $Z_0$, only. 
The past information (memory) $(Z_s)_{s\in (0,t)}$ is not necessary.
\end{rem}
\begin{proof}
By applying the Bayes rule and the relations (2.1) and (2.2), we have
\[
\pi_t(f):={\mathbb E}_y \left[ f(X)| {\mathcal F}^Z_t\right]
=\frac{\tilde{\mathbb E}_y \left[ \Lambda_t^{(y,T,X)} f(X)
\bigm| {\mathcal F}^Z_t\right]}
{\tilde{\mathbb E}_y \left[ \Lambda_t^{(y,T,X)} 
\bigm| {\mathcal F}^Z_t\right]},
\]
where we use the notation $\tilde{\mathbb E}_y[\cdot]$
for the expectation with respect to the probability measure
$d\tilde{\mathbb P}^1_{y} \otimes d{\mathbb P}_2$.
Observing that ${\mathcal F}^Z_t$ and $X$
are independent under $d\tilde{\mathbb P}^1_{y} \otimes d{\mathbb P}_2$, 
we deduce that 
\begin{equation}\label{CE}
\pi_t(f)=
\frac{\displaystyle 
\int_E \frac{\tilde{p}_{T-t}(Z_t, z)}{\tilde{p}_T(y,z)} f(z) \nu (dz)}
{\displaystyle \int_E \frac{\tilde{p}_{T-t}(Z_t, z')}
{\tilde{p}_T(y,z')} \nu (dz')}.
\end{equation}
\end{proof}
We note that once the initial value $Z_0$ is fixed, 
the process $(Z_t)_{t\in [0,T)}$ is ``Markovian'' in the following sense:
\begin{prop}\label{RMB-transdens}
For $0\le s<t<T$, (1.5) holds where
\[
P_{s,t}(x,dy|z_0)
=q(s,x; t,y|z_0) m(dy)
\]
and
\[
q(s,x; t,y|z_0)
:=\frac{\displaystyle
\tilde{p}_{t-s}(x,y) 
\left\{
\int_E \frac{\tilde{p}_{T-t}(y,z)}{\tilde{p}_T(z_0,z)}\nu(dz)
\right\}
}
{\displaystyle
\int_E \frac{\tilde{p}_{T-s}(x,z')}{\tilde{p}_T(z_0,z')}\nu(dz')
}.
\]
For a fixed $z_0$, the transition density $q(s,x;t,y|z_0)$ satisfies 
the Chapman-Kolmogorov identity:
\[
q(s,x;u,z|z_0)=\int_E q(s,x;t,y|z_0)q(t,y;u,z|z_0) m(dy)
\]
for
$0< s< t< u<T$ and $x,y,z\in E$.
\end{prop}

\begin{rem}
Let $m(dx)=dx$ and $\nu(dx)=g(x)dx$. 
Then, one can check that (a)
\[
\lim_{t\uparrow T} 
\int_E \frac{\tilde{p}_{T-t}(y,z)}{\tilde{p}_T(Z_0,z)}\nu(dz)
=\int_E \frac{\delta_y(z)}{\tilde{p}_T(Z_0,z)}g(z)dz
=\frac{g(y)}{\tilde{p}_T(Z_0,y)}
\]
where $\delta_y(\cdot)$ is the Dirac delta function
with point mass at $y$, 
and (b)
\[
\lim_{t\uparrow T} P_{s,t}(Z_s,dx|Z_0)
=\frac{\displaystyle
\frac{\tilde{p}_{T-s}(Z_s,y)}{\tilde{p}_T(Z_0,y)} \nu(dy)}
{\displaystyle
\int_E \frac{\tilde{p}_{T-s}(Z_s,z')}{\tilde{p}_T(Z_0,z')}\nu(dz')}
=\pi_s(dx).
\]
\end{rem}

\begin{proof}
The proof is similar to the one for Proposition 2.1.
We see that, for $0\le s\le t<T$, 
\[
{\mathbb E}_{z_0} \left[ f(Z_t)| {\mathcal F}^Z_s\right]
=\frac{\tilde{\mathbb E}_{z_0} \left[ \Lambda_t^{(y,T,X)} f(Z_t)
\bigm| {\mathcal F}^Z_s\right]}
{\tilde{\mathbb E}_{z_0} \left[ \Lambda_t^{(y,T,X)} 
\bigm| {\mathcal F}^Z_s\right]},
\]
where we use the Bayes rule, (2.1) and (2.2).
Recalling that ${\mathcal F}^Z_t$ and $X$
are independent under $d\tilde{\mathbb P}^1_{z_0} \otimes d{\mathbb P}_2$,  
the numerator of the right-hand side of the above equation is equal to
\begin{align*}
\int_{E} \tilde{p}_{t-s}(z_0,y) 
\left\{
\int_{E} 
\frac{\tilde{p}_{T-t}(y,z)}
{\tilde{p}_T(z_0,z)} \nu(dz) 
\right\} f(y) m(dy).
\end{align*} 
By computing the denominator in a similar way, we obtain the expression 
for $P_{s,t}(x,dy|z_0)$.
The Chapman-Kolmogorov identity is seen from 
\[
{\mathbb E}[f(Z_u)|{\mathcal F}^Z_s]
={\mathbb E}\left[ {\mathbb E}\left[ f(Z_u)
|{\mathcal F}^Z_t\right]
\bigm|{\mathcal F}^Z_s\right],
\quad
0\le s\le t\le u<T, 
\] 
which is the tower property of conditional expectation.
\end{proof}

Next, we are interested in 
describing the $({\mathbb P},{\mathcal F}^Z_t)$-dynamics of 
$\left(\pi_t(f)\right)_{t\in [0,T]}$. 
To this end, we construct the RMB $(Z_t)_{t\in [0,T]}$, 
starting with the underlying (unconditioned) Markovian SDE,
\begin{equation}\label{Y-diffusion}
 dY_t=b(Y_t)dt + \sigma(Y_t)d\tilde{W}_t,
\quad Y_0=y
\end{equation}
on $(\Omega^1_T, {\mathcal F}^1_T, \tilde{\mathbb P}^1_y,
({\mathcal F}^1_t)_{t\in [0,T]})$.
Here,  
$b:{\mathbb R}^n \to {\mathbb R}^n$
and $\sigma:{\mathbb R}^n \to {\mathbb R}^{n\times d}$
are given functions.
Further, $\tilde{W}:=\bigl(\tilde{W}^1,\dots,\tilde{W}^{d}\bigr)^\top$,
$\tilde{W}^i:=\bigl(\tilde{W}^i_t\bigr)_{t\in [0,T]}$ 
is a standard $d$-dimensional 
$(\tilde{\mathbb P}^1_y, {\mathcal F}^1_t)$-Brownian motion. Although assuming the existence (in the sense of its law) of a unique weak solution to (\ref{Y-diffusion}) is sufficient for constructing an RMB, we assume that there exists a strong, pathwise unique solution to (\ref{Y-diffusion}) in order to consider the filtering problem with respect to ${\mathcal F}^Z_t$.
Moreover, we assume that there exist associated transition densities
$\tilde{p}_t(x,y)>0$
with respect to the Lebesgue measure 
for $t\in (0,T]$, $x,y\in E$, 
which are sufficiently smooth with respect to $(t,x)$.
We recall that $\tilde{p}_{(\cdot)}(\cdot,y)$
satisfies the Kolmogorov backward equation 
\begin{equation}
\begin{split}
\left( -\partial_t + {\mathcal L}\right)\tilde{p}_t(x,y)=&0 
\quad\text{for $(t,x)\in (0,T] \times E$,} \\
\tilde{p}_0(x,y)=&\delta_y(x),
\end{split}
\end{equation}
where 
\[
{\mathcal L}(f)(x):= b(x)^\top \nabla f(x) + 
\frac{1}{2} {\rm tr}\left(
 (\sigma\sigma^\top)(x)\nabla\nabla f(x)\right) 
\]
with
$\nabla f:=\left( \partial_{x_1}f,\dots,\partial_{x_n} f\right)^\top$
and $\nabla\nabla f:=\left( \partial_{x_i x_j}f\right)_{1\le i,j\le n}$
being the infinitesimal generator of the Markov process $Y$.
In this case, 
by using It\^o's formula, we see that (2.2) satisfies
\[
d\Lambda^{(y,T,x)}_t =\Lambda^{(y,T,x)}_{t}
\left[
\nabla \ln \tilde{p}_{T-t}(Y_t,x)^\top \sigma(Y_t)d\tilde{W}_t
\right],
\]
where we denote $\nabla \ln \tilde{p}_t(y,x):=
\left( \partial_{y_1} \ln \tilde{p}_t(y,x),\dots,
\partial_{y_n} \ln\tilde{p}_t(y,x) \right)^\top$.
Hence, by Girsanov's theorem, we deduce that
\[
 W_t:= \tilde{W}_t -
\int_0^t \sigma(Z_s)^\top 
\nabla \ln \tilde{p}_{T-s}(Z_s,X) ds,
\quad t\in [0,T)
\]
is a $({\mathbb P}_y, {\mathcal F}_t)$-Brownian motion, 
and the SDE, defined on 
$(\Omega,{\mathcal F}_{T-}, {\mathbb P}_{y}, 
({\mathcal F}_t)_{t\in [0,T)})$ 
and satisfied by $Z:=(Z_t)_{t\in [0,T)}$, can be written down as
\begin{equation}
dZ_t 
=\left\{ b(Z_t)+ 
(\sigma\sigma^\top)(Z_t) \nabla \ln \tilde{p}_{T-t}(Z_t,X) 
\right\} dt 
+\sigma(Z_t)dW_t
\quad (Z_0=y). 
\end{equation}
We deduce that
$Z_{T-}:=\lim_{t\uparrow T} Z_t=X$,
${\mathbb P}_y$-almost surely. 
For the calculation of the 
$({\mathbb P}_{y}, {\mathcal F}^Z_t)$-dynamics of 
$Z$ and $\pi_t(f)$, 
we prepare the following: (a) We rewrite (2.3) as
\begin{equation}
\pi_t(f)=\frac{\rho_t(f)}{\rho_t(1)},
\end{equation}
where 
\[
 \rho_t(f):= \int_{E} f(z)\rho_t(dz)
\quad\text{and}\quad
 \rho_t(dz):=\frac{\tilde{p}_{T-t}(Z_t,z)}{\tilde{p}_T(y,z)} \nu(dz).
\]
(b) We define
\begin{equation}
R(t,z;f):= \int_{E} f(z') 
\frac{\tilde{p}_{T-t}(z,z')}{\tilde{p}_T(y,z')} \nu(dz') 
\quad\text{and}\quad
\Pi(t,z;f):=\frac{R(t,z;f)}{R(t,z;1)},
\end{equation}
where we fix $Z_0=y$, 
recalling the relations
\[
\rho_t(f)=R(t,Z_t;f)
\quad\text{and}\quad
\pi_t(f)=\Pi(t,Z_t;f)
=\frac{R(t,Z_t;f)}{R(t,Z_t;1)}.
\]
With this at hand, we obtain the subsequent result.
\begin{prop}\label{prop3}
Let $f:E\to {\mathbb R}$ be Borel-measurable so that
$\int_E |f(x)|\nu(dx)<\infty$.
Assume that $(\sigma\sigma^\top)>0$
and that
\begin{equation}
\int_0^t \int_{E} 
\left|
\sigma^\top(Z_s)
\nabla \tilde{p}_{T-s}(Z_s,x)
\right|^2 
\frac{f(x)}{\tilde{p}_T(Z_0,x)} 
\nu(dx) ds<\infty
\end{equation}
for any $t\in [0,T)$.
For $(s,x)\in (0,T)\times E$, let
\[
 \ell_s(x):=\ln \tilde{p}_{T-s}(Z_s,x).
\]
Then, the following statements hold:

\noindent{\rm (1)} 
On the filtered probability space
$(\Omega,{\mathcal F}_{T-}, {\mathbb P}_{y}, 
({\mathcal F}^Z_t)_{t\in [0,T)})$,  
the RMB $Z:=(Z_t)_{t\in [0,T)}$ satisfies
\begin{equation}
dZ_t 
=\left\{ b(Z_t)+ 
(\sigma\sigma^\top)(Z_t) 
\pi_t\left( \nabla \ell_t\right) 
\right\} dt
+(\sigma\sigma^\top)^{1/2}(Z_t)dU_t,
\quad 
(Z_0=y), 
\end{equation}
and $Z_{T-}:=\lim_{t\uparrow T} Z_t=X$
${\mathbb P}_y$-almost surely. Here, 
\begin{align}
 U_t:=&\int_0^t (\sigma\sigma^\top)^{-1/2} (Z_s)
\left\{ 
dZ_s -b(Z_s)ds
-(\sigma\sigma^\top)(Z_s)\pi_s (\nabla\ell_s) ds \right\} 
\nonumber \\
=&\int_0^t (\sigma\sigma^\top)^{-1/2} (Z_s)
\left\{ \sigma(Z_s) d\tilde{W}_s
-(\sigma\sigma^\top)(Z_s)\pi_s (\nabla\ell_s) ds \right\}
\end{align}
is an $n$-dimensional  
$({\mathbb P}_{y}, {\mathcal F}^Z_t)$-Brownian motion. 

\noindent{\rm (2)} 
For $t\in [0,T)$, it holds that
\begin{equation}
 \pi_t(f)=\pi_0(f)
+ \int_0^t 
\left\{
\pi_s\left(f\nabla\ell_s\right)-\pi_s(f)\pi_s\left( \nabla\ell_s\right)
\right\}^\top 
(\sigma\sigma^\top)^{1/2}(Z_s)dU_s.
\end{equation}
Moreover, we have 
\begin{equation}
 \nabla \Pi(s,z;f)=
\Pi(s,z; f\nabla \ell_s)
-\Pi(s,z;f)\Pi(s,z;\nabla \ell_s), 
\end{equation}
where $\nabla\Pi(s,z;f)^\top:=
\left( \partial_{z_1}\Pi,\dots,\partial_{z_n}\Pi\right)(s,z;f)$, 
and hence
\begin{equation}
 \pi_t(f)=\pi_0(f)
+ \int_0^t (\nabla \Pi)(s,Z_s;f)^\top 
(\sigma\sigma^\top)^{1/2}(Z_s)dU_s
\end{equation}
follows.
\end{prop}

\begin{rem}
From a stochastic filtering viewpoint, 
(2.12)--(2.14) are recognised as the Kushner-Stratonovich equations
for $\left(\pi_t(f) \right)_{t\in [0,T)}$, 
which has the finite-dimensional Markovian ``information state'' 
process $Z$.
\end{rem}

\begin{rem}
Similar SDEs to (2.10) appear in 
Section 3 of Baudoin (2002) and 
Theorem 2.1 and 3.1 of \c{C}etin \& Danilova (2016). 
\end{rem}

\begin{proof}
(1) By the Kolmogorov backward equation (2.5) and It\^o's formula, 
we have
\[
\tilde{p}_{T-t}(Z_t,x)
=
\tilde{p}_{T}(Z_0,x) 
+\int_0^t \nabla \tilde{p}_{T-s}(Z_s,x)^\top 
\sigma(Z_s)d\tilde{W}_s.
\]
Then, by making use of this relation, we deduce that 
\begin{align*}
 \rho_t(f)
=& \int_{E} f(x)
\frac{\tilde{p}_{T-t}(Z_t,x)}{\tilde{p}_T(Z_0,x)} \nu(dx) \\
=& \rho_0(f)
+ \int_{E} \frac{f(x)}{\tilde{p}_T(Z_0,x)} 
\left[
\int_0^t \tilde{p}_{T-s}(Z_s,x)\nabla \ln \tilde{p}_{T-s}(Z_s,x)^\top 
\sigma(Z_s)d\tilde{W}_s
\right] \nu(dx) \\
=& \rho_0(f)
+\int_0^t 
\left[
\int_E f(x)
\nabla \ln \tilde{p}_{T-s}(Z_s,x)^\top  
\frac{\tilde{p}_{T-s}(Z_s,x)}{\tilde{p}_T(Z_0,x)}
\nu(dx) 
\right]
\sigma(Z_s)d\tilde{W}_s \\
=& \rho_0(f)
+\int_0^t \rho_s(f\nabla \ell_s^\top) \sigma(Z_s) d\tilde{W}_s .
\end{align*}
Here, the stochastic version of Fubini's theorem is employed 
while recalling the condition (2.9),
we refer to Theorem 4.65 in Protter (2005). 
So, we see that
\[
\rho_t(1)=
 \tilde{\mathbb E}_y\left[ \Lambda^{(y,T,X)}_t \Bigm| {\mathcal F}^Z_t\right]
\]
satisfies
\[
d\rho_t(1)=
\rho_t(\nabla \ell_t^\top) \sigma(Z_t) d\tilde{W}_t 
=\rho_{t}(1)
\left\{
\pi_t(\nabla \ell_t^\top) \sigma(Z_t) d\tilde{W}_t
\right\},
\]
where we recall (2.7). 
Using the above conditional 
$\left( \tilde{\mathbb P}_y, {\mathcal F}^Z_t\right)$-martingale 
density representation,
we apply Girsanov's theorem to see that
\begin{align*}
 &\int_0^t \sigma(Z_u) 
\left\{d\tilde{W}_u-
\sigma(Z_u)^\top \pi_u(\nabla \ell_u) du
\right\} \\
=&
Z_t -Z_0-\int_0^t b(Z_u)du
-\int_0^t (\sigma\sigma^\top)(Z_t) \pi_t(\nabla \ell_t) du,
\quad t\in [0,T) 
\end{align*}
is a $\left( {\mathbb P}_y, {\mathcal F}^Z_t\right)$
local-martingale. Hence, we deduce that
$(U_t)_{t\in [0,T)}$ given by (2.11) is 
a $\left( {\mathbb P}_y, {\mathcal F}^Z_t\right)$-Brownian motion
by using L\'evy's characterization.
(2.10) is seen by combining (2.4) and (2.11).

\medskip

\noindent{(2)}
By using It\^o's formula, (2.7) and (2.11), we see that
\begin{align*}
d\pi_t(f)=&\frac{d\rho_t(f)}{\rho_t(1)}
-\frac{\rho_t(f) d\rho_t(1)}{\rho_t(1)^2}
-\frac{d\langle \rho(f), \rho(1)\rangle_t}{\rho_t(1)^2}
+\frac{\rho(f) d\langle \rho(1)\rangle_t}{\rho_t(1)^3} \\
=&\frac{\rho_t(f\nabla \ell_t^\top)\sigma(Z_t)d\tilde{W}_t}{\rho_t(1)}
-\frac{\rho_t(f)\rho_t(\nabla\ell_t^\top)\sigma(Z_t)d\tilde{W}_t}{\rho_t(1)^2} \\
&-\frac{\rho_t(f\nabla\ell_t^\top) (\sigma\sigma^\top)(Z_t) 
\rho_t(\nabla\ell_t)}{\rho_t(1)^2}dt
+\frac{\rho_t(f)\rho_t(\nabla\ell_t^\top) (\sigma\sigma^\top)(Z_t) 
\rho_t(\nabla\ell_t)}{\rho_t(1)^3}dt \\
=& 
\left\{
\pi_t(f\nabla\ell_t^\top)
-\pi_t(f)\pi_t(\nabla\ell_t^\top)
\right\}
\bigl(\sigma\sigma^\top\bigr)^{1/2}(Z_t) dU_t,
\end{align*}
from which (2.12) follows. Equation
(2.13) can be verified directly, and hence (2.14) follows.
\end{proof}

\begin{rem}
If we admit $(t,z)\mapsto \Pi(t,z;f)$ is in $C^{1,2}([0,T)\times E)$,
then (2.14) is obtained from It\^o's formula.
\end{rem}

The following is a direct interpretation of Proposition 2.3
in a financial context.
\begin{cor}
On the filtered probability space
$(\Omega,{\mathcal F}_{T-},{\mathbb P}_y,({\mathcal F}^Z_t)_{t\in [0,T)})$,
the asset price process $(S_t)_{t\in [0,T]}$ given by (1.4)
has the dynamics,
\[
 dS_t =r S_t dt+e^{-r(T-t)}\Sigma_t(\sigma\sigma^\top)^{1/2}(Z_t) dU_t,
\quad t\in [0,T).
\]
Here, $(U_t)_{t\in [0,T)}$ is a $({\mathbb P}_y, {\mathcal F}^Z_t)$-Brownian
motion given by (2.11), 
\[
\Sigma_t:=\nabla \Pi(t,Z_t;f),
\]
where we use (2.8), and 
the dynamics of $(Z_t)_{t\in [0,T)}$ is given by (2.10).
\end{cor}

\begin{rem}
As can be seen in Proposition 2.3, 
the volatility term $\Sigma_t$ 
can be rewritten as
\begin{equation}
 \Sigma_t=
\left\{
\pi_t\left( f \nabla \ell_t\right)
-\pi_t(f)\pi_t\left(\nabla \ell_t\right)
\right\}^\top.
\end{equation}
In the case of the linear Gaussian diffusion (3.1) treated in the next section, 
(2.15) is simplified to 
\begin{equation}
 \Sigma_t=
{\rm Cov}_y\left[ f(X),X \bigm| {\mathcal F}^Z_t\right]^\top
V_{T-t}^{-1} e^{(T-t)K},
\end{equation}
where we use (3.3) and (3.12), and recall that
\begin{equation}
{\rm Cov}_y\left[ f(X), X \bigm| {\mathcal F}_t^Z\right] 
:= {\mathbb E}_y\left[ f(X)X\bigm| {\mathcal F}^Z_t\right]
-{\mathbb E}_y\left[ f(X)\bigm| {\mathcal F}^Z_t\right]
{\mathbb E}_y\left[ X | {\mathcal F}^Z_t\right].
\end{equation}
We leave it to later to give an economic and financial interpretation 
to the representation (2.16) of the volatility process. In Section 4.2, Proposition 4.1, we apply (2.16) to compute the price sensitivities 
of a financial security linked to greenhouse gas emissions
with respect to model parameters.
\end{rem}
We end this section by giving a simple example of an RMB, 
which is introduced in Brody et al. (2008).
\begin{ex}[Randomized Brownian Bridge]
In (2.4), let $n=1$ and let $b\equiv 0$ and $\sigma\equiv 1$. That is, 
we start with
\[
 Y_t =y+ W_t,
\quad t\ge 0
\]
on the state space $E:={\mathbb R}$, where its transition density 
with respect to the Lebesgue measure $m$ is written as
\[
 \tilde{p}_t(y,z)=\frac{1}{\sqrt{2\pi t}} 
e^{-\frac{|y-z|^2}{2t}}.
\]
Then, recalling (2.6) and
\[
 \nabla \ell_t(x)=\frac{1}{T-t}\left( x-Z_t\right),
\]
the SDE for $Z:=Y^{(y,T,X)}$ on
$(\Omega,{\mathcal F}_{T-},{\mathbb P}_y,({\mathcal F}_t)_{t\in [0,T)})$
is given by 
\[
 dZ_t =\frac{X-Z_t}{T-t} dt +dW_t,
\quad Z_0=y,
\]
whose solution may be written as
\[
 Z_t = \frac{t}{T} X
+\frac{T-t}{T}\left( y + \int_0^t \frac{T}{T-s}dW_s\right).
\]
We see that, 
in this case, the hidden variable $X$ and the noise part 
are separated in an additive way. 
The conditional probability 
$\pi_t(dz)=\pi_t(dz|Z_t,Z_0)$ is represented as a functional of $(Z_t,Z_0)$
as shown in Proposition 2.1, 
and the dynamics of $Z$ 
on $(\Omega,{\mathcal F}_{T-},{\mathbb P}_y,({\mathcal F}^Z_t)_{t\in [0,T)})$
is given by
\[
dZ_t = \frac{\pi_t(I)-Z_t}{T-t}dt +dU_t,
\]
where we recall (2.10), use 
$I(x):=x$ and the ${\mathcal F}^Z_t$-Brownian motion 
$(U_t)_{t\in [0,T)}$ defined by (2.11).
\end{ex}

\section{Skew-normal randomised diffusion bridge}
In this section, we introduce a new example of an RMB, where 
the conditional probability $\pi_t(dz)$ given in Proposition 2.1 
and the transition density $q(s,x;t,y|z_0)$ given in Proposition 2.2 
have explicit representations. 
\

As the underlying Markov process, we employ 
the linear Gaussian diffusion
\begin{equation}\label{Ydiff}
 dY_t = (k+ KY_t)dt + \Sigma dW_t,
\quad Y_0\in {\mathbb R}^n,
\end{equation}
where $k\in {\mathbb R}^n$, 
$K$ and $\Sigma\in {\mathbb R}^{n\times n}$ so that 
$\Sigma\Sigma^\top >0$, 
and where $(W_t)_{t\ge 0}$
is an $n$-dimensional $({\mathcal F}_t)$-Brownian motion.
In this case, the transition density $\tilde{p}_t(y,z)$ that satisfies
\[
 {\mathbb P}\left( Y_t\in dz | Y_s=y\right)
= \tilde{p}_{t-s}(y,z)dz
\]
for $0\le s\le t$, is given by
\begin{equation}
\tilde{p}_{t} (y,z)=\phi_n\left( z- \mu_t(y), V_t\right)
\end{equation}
where 
\[
 \phi_n \left( z; V_t\right)
:=\frac{1}{(2\pi)^{\frac{n}{2}}\sqrt{\det(V_t)}}
\exp\left(-\frac{1}{2}z^\top V_t^{-1} z\right)
\]
is the probability density function
of the $n$-dimensional normal distribution 
with the mean vector $\mu_t\in {\mathbb R}^n$
and the covariance matrix $V_t\in {\mathbb S}^n_{++}
:=\left\{ S\in {\mathbb S}^n| S>0\right\}$, 
and
\begin{equation}\label{meanvar}
\begin{split}
\mu_t(y):=& e^{tK} y
+ 
\left(\int_0^t e^{sK} ds\right)
k, \\
V_t:=& \int_0^t e^{sK} \Sigma\Sigma^\top e^{sK^\top} ds.
\end{split}
\end{equation}
Next, we set the law of the random variable $X$ as
the generalised multivariate skew normal distribution
${\rm GMSN}(a,b,A,B,C)$. That is:
\begin{equation}
\nu(dx):={\mathbb P}(X\in dx)
:= f_{\rm GMSN}(x; a,b,A,B,C) dx,
\end{equation}
where, for $a\in {\mathbb R}^m$,
$b\in {\mathbb R}^n$, 
$A\in {\mathbb S}^m_{++}$, 
$B\in {\mathbb S}^n_{++}$, 
$C\in {\mathbb R}^{m\times n}$, 
$m\in {\mathbb N}$, $n\in {\mathbb N}$, and
\[
 \Phi_m(x; A)
:= \int_{\prod_{i=1}^m(-\infty,x_i]}\phi_m(y;A)dy,
\]
we have 
\begin{equation}\label{fGMSN}
f_{\rm GMSN}(x; a,b,A,B,C) 
:= 
\frac{1}{\displaystyle \Phi_m\left( Cb- a,A+CBC^\top\right)}
\phi_n(x-b,B)
\Phi_m\left( Cx- a, A\right).
\end{equation} 
\begin{rem}
The distribution ${\rm GMSN}(a,b,A,B,C)$ 
is introduced in Section 5 of Gupta et al. (2004). 
In the case
\[
\begin{pmatrix}
 X_1 \\ X_2
\end{pmatrix}
\sim
{\mathcal N}
\left(
\begin{pmatrix}
a \\ b
\end{pmatrix},
\begin{pmatrix}
A+ C B C^\top& -CB \\
-B C^\top& B
\end{pmatrix}
\right), 
\] 
we see that $X_2|\{ X_1\le Cb\} \sim {\rm GMSN}(a,b,A,B,C)$.
Here, noting the relation
\[
{\rm GMSN}(a,b,A,B,0)={\mathcal N}(b,B),  
\]
we see that the parameter $C$ controls the skewness of the GMSN distribution.
\end{rem}

\begin{rem}[Multivariate Skew Normal Distribution]
Indeed, Gupta et al. (2004) mainly focus on 
the multivariate skew normal distribution (MSN),  
where $n=m$ and
\begin{align*}
f^{(1)}_{\rm MSN}\left( x; b,B,C\right)
:=&f_{\rm GMSN}(x;Cb,b,I_n,B,C) \\
=&\frac{1}{\displaystyle \Phi_n\bigl( 0; I_n+C B C^\top\bigr)}
\phi_n(x-b;B)\Phi_n\left( C(x-b);I_n \right).
\end{align*}
Azzalini and Dalla Valle (1996)
and Azzalini and Capitanio (1999, 2013)
focus on yet another MSN distribution,
where $m=1$, $b\in {\mathbb R}^n$, $c\in {\mathbb R}^n$, and 
\begin{align*}
 f^{(2)}_{\rm MSN}(x;b,c):=&f_{\rm GMSN}(x;c^\top b,b,1,I_n,c^\top) \\
=&2\phi_n(x-b; I_n)\Phi_1\left( c^\top (x-b);1\right).
\end{align*}
Developments of these MSN families are 
motivated by modelling the skewness of probabilitistic models. 
For related studies, we refer to Azzalini and Capitanio (2013), for example.
\end{rem}

\begin{rem}[Univariate Skew Normal Distribution]
The univariate skew normal probability density is given by 
letting $m=n=1$ and by setting  
\begin{align*}
 f_{\rm USN}(x;\mu,\sigma,\alpha)
:=&f_{\rm GMSN}\left( x; \frac{\alpha\mu}{\sigma},\mu,1,\sigma,
\frac{\alpha}{\sigma}\right) \\
=& \frac{2}{\sigma}
\phi_1\left( \frac{y-\mu}{\sigma};1\right)
\Phi_1\left( \alpha \left( \frac{y-\mu}{\sigma}\right);1\right).
\end{align*}
The parameters $\mu\in {\mathbb R}$, $\sigma>0$, and 
$\alpha\in {\mathbb R}$
are called
the location parameter, 
the scale parameter, 
and the shape parameter, respectively. 
For various properties and related statistical methodologies of USN, 
we refer to, e.g., Chapter 1-4 of Azzalini and Capitanio (2013).
\end{rem}

\begin{prop}
For the RMB specified by (3.1)-(3.5), 
the conditional probability $\pi_t(dz)$ given in Proposition 2.1
is conditionally GMSN. That is, for $t\in [0,T)$,
\[
 \pi_t(dz)
=f_{\rm GMSN}\left( z; 
a,\tilde{b}(t,Z_0,Z_t),A,\tilde{B}_t,C \right)dz,
\]
where
\begin{align*}
\tilde{b}(t,Z_0,Z_t):=& 
\left( V_{T-t}^{-1}-V_T^{-1}+B^{-1}\right)^{-1}
\left\{ V_{T-t}^{-1}\mu_{T-t}(Z_t)
-V_T^{-1}\mu_{T}(Z_0)+B^{-1}b\right\}, \\
\tilde{B}_t:=& \left( V_{T-t}^{-1}-V_T^{-1}+B^{-1}\right)^{-1}.
\end{align*}
\end{prop}
The following lemma is useful for proving the proposition.
\begin{lem}
Let $P\in {\mathbb S}^n_{+}:=\left\{ S\in {\mathbb S}^n| S\ge
 0\right\}$ and $p\in {\mathbb R}^n$.
For the GMSN density given by (3.5), it holds that
\begin{align*}
&\exp\left( -\frac{1}{2}y^\top P y + p^\top y\right)
f_{\rm GMSN}(y; a,b,A,B,C) \\
&=
f_{\rm GMSN}\left(y; a,\tilde{b},{A},\tilde{B},C\right) \\
&\quad\times\frac{1}{\sqrt{\det(I+B P)}}
\frac{\displaystyle 
\Phi_m\left( C\tilde{b}-a
; A+C\tilde{B} C^\top\right)}
{\displaystyle \Phi_m\left( Cb-a; A+CB C^\top\right)}
\exp\left\{
\frac{1}{2}\left(
\tilde{b}^\top\tilde{B}^{-1}\tilde{b}
-b^\top B^{-1}b\right)\right\},
\end{align*}
where
$\tilde{b}:=(P+B^{-1})^{-1}(p+B^{-1}b)$ and  $\tilde{B}:=(P+B^{-1})^{-1}$.
\end{lem}

\begin{proof}
We see that
\begin{align}
&\exp\left( -\frac{1}{2}y^\top P y + p^\top y\right)
\phi_n(y-b;B) 
\nonumber \\
&= 
\frac{1}{(2\pi)^{\frac{n}{2}}\sqrt{\det(B)}} 
\exp\left\{
-\frac{1}{2}
(y-\tilde{b})^\top \tilde{B}^{-1} (y-\tilde{b})
+\frac{1}{2}\left(
\tilde{b}^\top\tilde{B}^{-1}\tilde{b}
-b^\top B^{-1}b\right)
\right\} \nonumber \\
&=\phi_n\left( y-\tilde{b}; \tilde{B}\right) 
\frac{1}{\sqrt{\det(I+B P)}}
\exp\left\{\frac{1}{2}\left(
\tilde{b}^\top\tilde{B}^{-1}\tilde{b}
-b^\top B^{-1}b\right)\right\}.
\end{align}
Furthermore,
\begin{align}
\phi_n\left( y-\tilde{b};\tilde{B}\right) 
\Phi_m\left( Cy-a; A\right)
=
\Phi_m\left( C\tilde{b}-a; A+C\tilde{B} C^\top\right)
f_{\rm GMSN}
\left( y; a, \tilde{b}, A,\tilde{B},C\right).
\end{align}
The proof is completed by combining (3.6) and (3.7).
\end{proof}
Next we prove Proposition 3.1:
\begin{proof}[Proof]
We observe that
\begin{align*}
\frac{p_{T-t}(y,z)}{p_T(y_0,z)}
&=\sqrt{\frac{\det V_T}{\det V_{T-t}}} 
\exp\left[
-\frac{1}{2}
\left\{ z-\mu_{T-t}(y)\right\}^\top 
V_{T-t}^{-1}\left\{ z-\mu_{T-t}(y)\right\}
\right. \\
&\hspace{4cm}\left. +\frac{1}{2}
\left\{ z-\mu_{T}(y_0)\right\}^\top 
V_{T}^{-1}\left\{ z-\mu_{T}(y_0)\right\}
\right] \\
&=\sqrt{\frac{\det V_T}{\det V_{T-t}}} 
\exp\left[
-\frac{1}{2}z^\top (V_{T-t}^{-1}-V_T^{-1}) z
+\left\{ V_{T-t}^{-1}\mu_{T-t}(y) 
-V_T^{-1} \mu_T(y_0)\right\}^\top z \right. \\
&\hspace{4cm}\left.
-\frac{1}{2}\mu_{T-t}(y)^\top V_{T-t}^{-1}\mu_{T-t}(y)
+\frac{1}{2}\mu_{T}(y_0)^\top V_{T}^{-1}\mu_{T}(y_0)
\right].
\end{align*} 
We then apply Lemma 3.1 to obtain
\begin{align}
J_{t,T}(z;y_0,y) 
:=&\frac{p_{T-t}(y,z)}{p_T(y_0,z)}\frac{d\nu}{dz}(z) 
\nonumber \\
=& 
f_{\rm GMSN}\left(z; a,\tilde{b}(t,y_0,y),{A},\tilde{B}_t,C\right) 
\nonumber \\
\times&\frac{1}{\sqrt{\det\left(I+{B}(V_{T-t}^{-1}-V_T^{-1})\right)}}
\frac{\displaystyle 
\Phi_m\left( C\tilde{b}(t,y_0,y)-a;A+C\tilde{B}_t C^\top\right)}
{\displaystyle \Phi_m\left( Cb-a; A+CB C^\top\right)} 
\nonumber \\
\times&
\sqrt{\frac{\det V_T}{\det V_{T-t}}} 
\exp\Biggl[
\frac{1}{2}\biggl\{
\left(\tilde{b}^\top\tilde{B}_t^{-1}\tilde{b}\right)(t,y_0,y)
-b^\top B^{-1}b 
\nonumber \\
&\qquad-\mu_{T-t}(y)^\top V_{T-t}^{-1}\mu_{T-t}(y)
+\mu_{T}(y_0)^\top V_{T}^{-1}\mu_{T}(y_0)
\biggr\} \Biggr].
\end{align}
Hence, the statement of Proposition 3.1 follows by recalling that
\[
\pi_t(dz)=
\frac{J_{t,T}(z;Z_0,Z_t)dz}
{\displaystyle\int_{{\mathbb R}^n}J_{t,T}(z;Z_0,Z_t)dz}.
\]
\end{proof} 
\begin{rem}
By letting $P=0$ in Lemma 3.1, we obtain the explicit expression 
of the moment generating function of a GMSN-distributed multivariate random variable:
\begin{align*}
M(p):=&
\int_{{\mathbb R}^n} e^{p^\top y}
f_{\rm GMSN}(y; a,b,A,B,C) dy \\
=&
\frac{\displaystyle 
\Phi_m\left( C(b+Bp)-a; A+CB C^\top\right)}
{\displaystyle \Phi_m\left( Cb-a; A+CB C^\top\right)}
\exp\left( \frac{1}{2}p^\top B p+b^\top p\right).
\end{align*}
From this we see for example that 
\begin{align}
\nabla M(0) 
=&\int_{{\mathbb R}^n} y
f_{\rm GMSN}(y; a,b,A,B,C) dy \nonumber \\
=&b+ B C^\top \left(\frac{\nabla\Phi_m}{\Phi_m}\right)
\left( Cb-a; A+CB C^\top\right), 
\end{align}
where 
$\nabla \Phi_m(x):=\left( \partial_{x_1} \Phi_m(x),
\dots,\partial_{x_m}\Phi_m(x)\right)^\top$.
\end{rem}
We introduce the density
\begin{equation}
\nu(dx):={\mathbb P}(X\in dx)
= f_{\rm GMSN}\left(x; a,\mu_T(Z_0),A,V_T,C\right) dx,
\end{equation}
which is constructed from 
the marginal law of $Y_T$, ${\mathcal N}(\mu_T(Z_0), V_T)$,
by adding a skewness effect. 
Then, the expression of the conditional probability $\pi_t(dz)$
can be simplified as follows:
\begin{cor}\label{Cor}
For an RMB, consisting of (3.1)-(3.3) and (3.10), it holds that
\[
 \pi_t(dz)
=f_{\rm GMSN}\left( z; a,\mu_{T-t}(Z_t),A,V_{T-t},C \right)dz,
\quad t\in [0,T).
\]
In particular, the conditional probability $\pi_t(dz)$
is independent of $Z_0$.
\end{cor}
Moreover, we obtain the following:
\begin{prop}
For the RMB specified by (3.1)-(3.3) and (3.10), 
the following holds: 

\noindent{\rm (1)} The transition density
$q(s,x;t,y|Z_0)$, obtained in Proposition 2.2, is expressed by
\[
q(s,x;t,y|Z_0)
=q(s,x;t,y)
=f_{\rm GMSN}
\left( y; \tilde{a}_t,\mu_{t-s}(x), \tilde{A}_t,V_{t-s}, \tilde{C}_t\right),
\]
where
\begin{align*}
\tilde{a}_t:=&a-C \left( \int_0^{T-t} e^{u K} du \right) k, \\
\tilde{A}_t:=&A+C V_{T-t} C^\top, \\
\tilde{C}_t:=&C e^{(T-t)K},
\end{align*}
and is independent of $Z_0$.

\noindent{\rm (2)} The dynamical equation of the RMB $Z$ under
$(\Omega,{\mathcal F}_{T-}, {\mathbb P}_y, ({\mathcal F}_t^Z)_{t\in [0,T)})$
is given by
\begin{multline}
dZ_t = 
\left\{ 
k+ KZ_t
-(\Sigma\Sigma^\top)C^\top
\left( \frac{\nabla \Phi_m}{\Phi_m}\right)
\left( C\mu_{T-t}(Z_t)-a; A+C V_{T-t}C^\top\right)
\right\} dt \\
+(\Sigma\Sigma^\top)^{1/2}dU_t
\quad
\end{multline}
where $Z_0=y$ and $(U_t)_{t\in [0,T)}$ is the 
$({\mathcal F}^Z_t)$-Brownian motion defined by (2.11) 
for $b(z)=k+Kz$ and $\sigma(z)=\Sigma$.
\end{prop}

\begin{proof}
1. We apply Proposition 2.2 and Formula (3.8), where 
$b=\mu_T(z_0)$ and $B=V_{T}$ and observe that
\begin{align*}
q(s,x;t,y|z_0) 
=&\frac{\displaystyle p_{t-s}(x,y) \int_{{\mathbb R}^n} J_{t,T}(z; z_0,y)dz}
{\displaystyle \int_{{\mathbb R}^n} J_{s,T}(z; z_0,x)dz} \\
=&\phi_n\left( y-\mu_{t-s}(x);V_{t-s}\right)
\frac{\displaystyle \Phi_m\left( C\mu_{T-t}(y)-a; A+CV_{T-t}
 C^\top\right)}
{\displaystyle \Phi_m\left( C\mu_{T-s}(x)-a; A+CV_{T-s} C^\top\right)}.
\end{align*}
By recalling here that
\begin{eqnarray}
 \Phi_m\left( C\mu_{T-t}(y)-a; A+CV_{T-t} C^\top\right)
&=&\Phi_m\left( \tilde{C}_t y-\tilde{a}_t;\tilde{A}_t\right),\nn \\
 \Phi_m\left( C\mu_{T-s}(x)-a; A+CV_{T-s} C^\top\right)
&=& \Phi_m\left( \tilde{C}_t\mu_{t-s}(x)-\tilde{a}_t;
\tilde{A}_t+\tilde{C}_tV_{t-s} \tilde{C}_t^\top\right),\nn
\end{eqnarray}
the assertion follows. 

2. We compute the SDE (2.10) utilising the setting of this section.
We consider $b(t,z)=k+Kz$, $\sigma(z)=\Sigma$, and 
\begin{align}
 \ell_t(z)=& \log \tilde{p}_{T-t}(Z_t,z) \nonumber \\ 
=&-\frac{1}{2}\left( z-\mu_{T-t}(Z_t)\right)^\top V_{T-t}^{-1}
\left( z-\mu_{T-t}(Z_t)\right)
-\frac{1}{2} \log \det (V_{T-t}) -\frac{n}{2}\log (2\pi).
\end{align}
By Formula (3.9) and Corollary 3.1, we deduce that
\begin{align*}
\pi_t\left( \nabla \ell_t\right)
=&V_{T-t}^{-1}
\left\{\mu_{T-t}(Z_t)-\int_{{\mathbb R}^n} z \pi_t(dz)\right\} \\
=&
-C^\top
\frac{\nabla \Phi_m}{\Phi_m}\left( C\mu_{T-t}(Z_t)-a; 
A+C V_{T-t}C^\top \right), 
\end{align*}
and the assertion follows.
\end{proof}

In Figures 1 and 2, considering the setting in Proposition 3.2, 
that is, the skew-normal RMB specified by (3.1)-(3.3)
and (3.10), we draw its marginal densities 
\[
 \left\{ q(0,x_0; t,x)\right\}_{(t,x)\in [T_0,T]\times {\mathbb R}},
\]
where we set $x_0=0$, $T_0=0.01$, $T=1$. The parameter values of the RMB are set to
$n=1$, 
$k=0$, 
$K=-0.1$, 
$\Sigma=0.3$, 
$a=1$, $A=1$, and $C=20$.
\begin{figure}[H]
\begin{center}\label{Fig:1}
\includegraphics[angle=-90, scale=0.5]{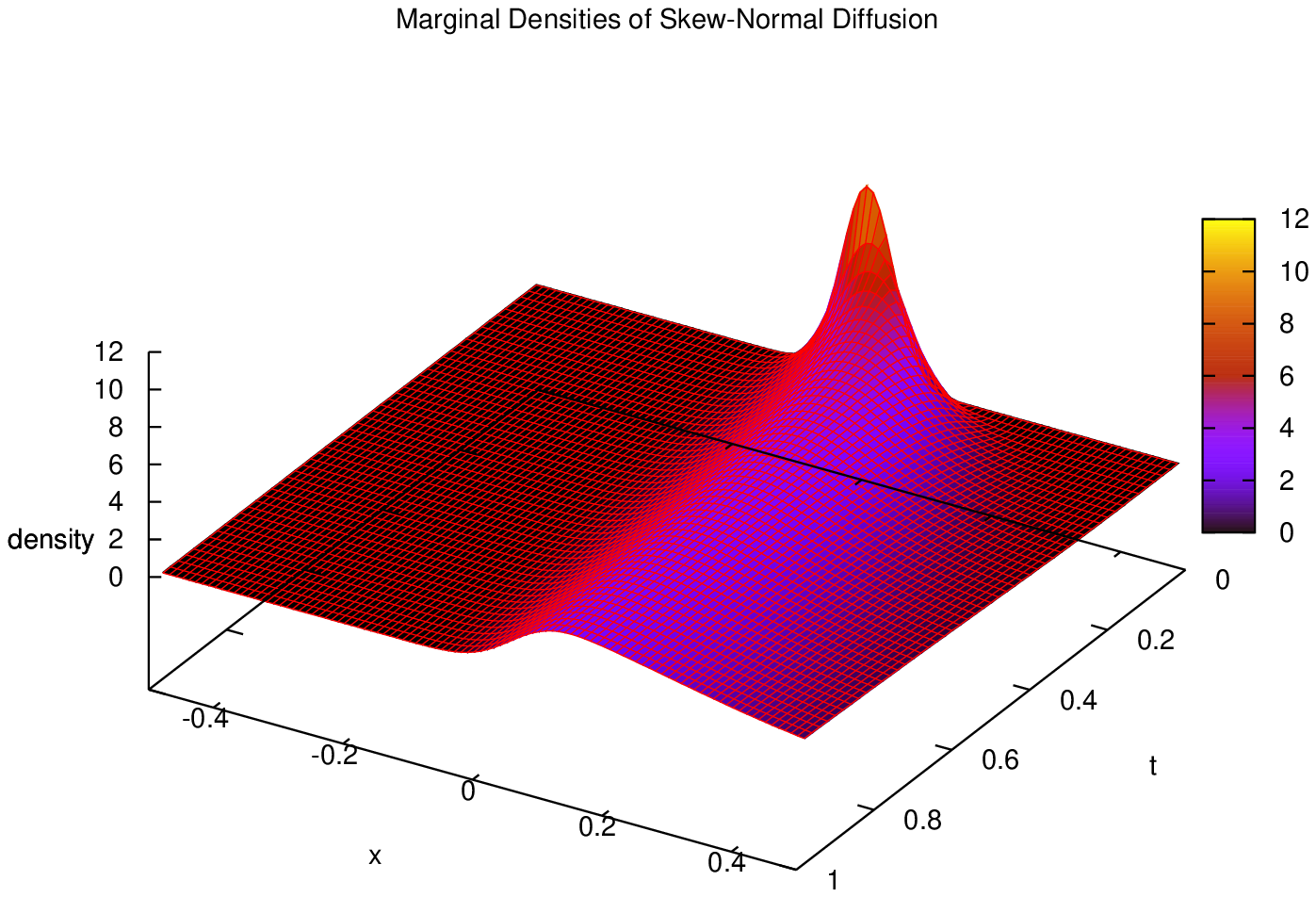}
\caption{Skew-normal densities $q(0,0;t,x)$ for $0.01\le t\le 1$.}
\end{center}
\end{figure}

\begin{figure}[H]
\begin{center}\label{Fig:2}
\includegraphics[angle=-90, scale=0.5]{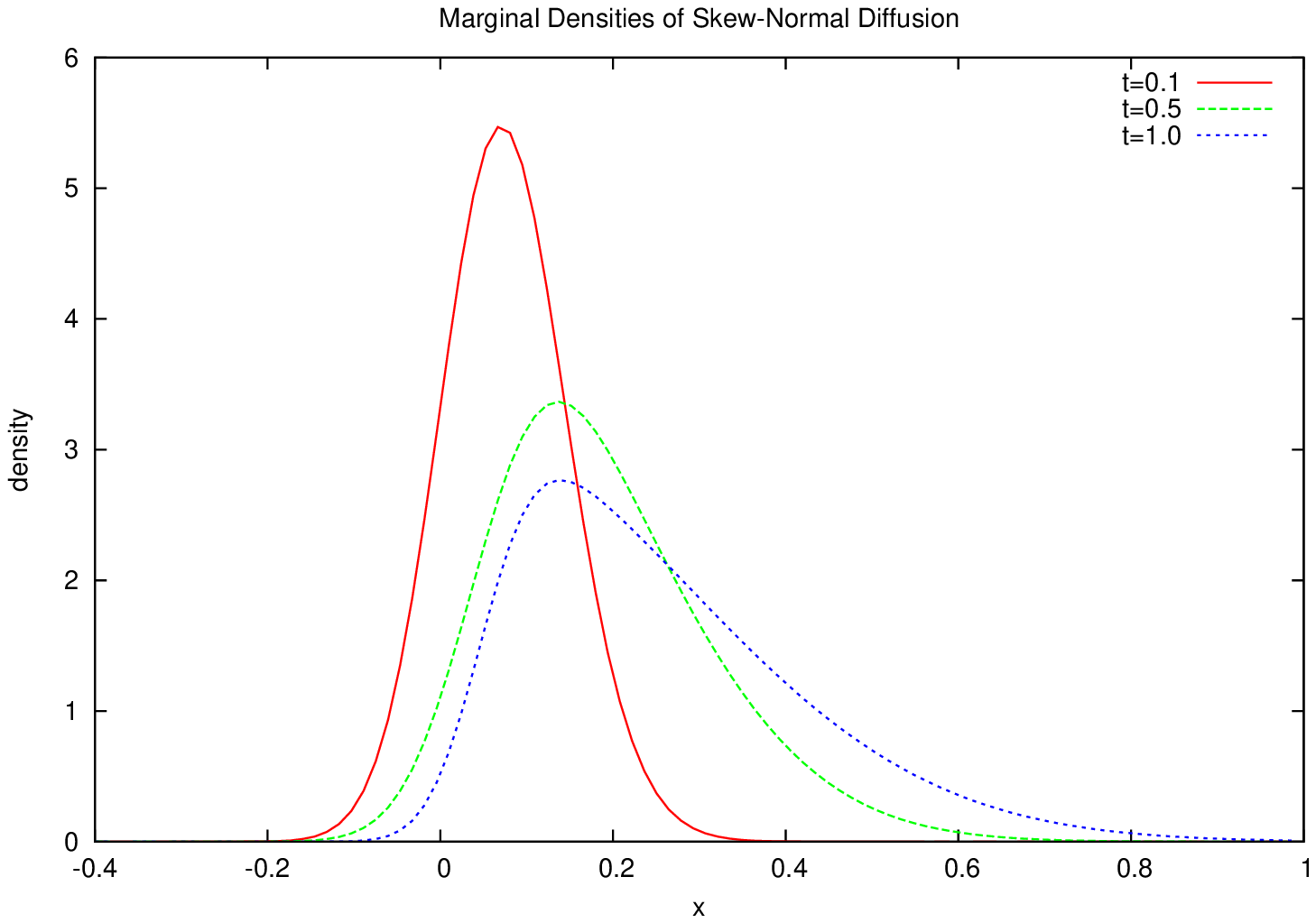}
\caption{Skew-normal densities $q(0,0;t,x)$ for $t\in\{0.1, 0.5, 1.0\}$.}
\end{center}
\end{figure}

\section{Applications to commodity pricing and securitization of greenhouse gas risk}
In this section we apply the randomised Markov bridges to construct asset price models following Formula (\ref{pricing_formula}). We consider skew-normal randomised diffusion bridges (SNRDB) to price commodity assets and for the securitisation of greenhouse gas (GHG) risk. The applications driven by SNRDB are just an example; the pricing models can be constructed with any RMBs, provided one agrees with the suitability of the implicit statistics for the envisaged application.  

It is customary in commodity pricing to assume a spot price process $(S_t)_{0\le t\le T}$ for the considered commodity and to model the price process $(F_{tT})_{0\le t\le T}$ of a $T$-delivery forward contract. The spot price process is $(\F_t)$-measurable and the delivery date $T$ is a fixed date specified in the forward contract. The price $F_{tT}$ at time $t$ of the forward contract with delivery date $T$ is given by
\begin{equation}\label{f-pp}
F_{tT}=\E\left[S_T\,\vert\,\F_t\,\right]
\end{equation}
where the expectation is taken under the risk-neutral measure. 
We assume  that the risk-free interest rate is constant,
which implies the distribution of $S_T$ under the $T$-forward measure 
is unchanged, and we can evaluate the price $F_{tT}$ 
under the forward measure $\PR^T$ by the relation
\begin{equation}\label{fwd-f-pp}
F_{tT}=
\E^{\PR^T}\left[S_T\,\vert\,\F_t\,\right]=\E\left[S_T\,\vert\,\F_t\,\right].
\end{equation}
Next, we construct a commodity pricing model driven by a randomised diffusion bridge with terminal GMSN marginal.
\subsection{Commodity pricing}
We take the RMB $(Z_t)_{t\in [0,T]}$, specified by (3.1)-(3.3) and (3.10), as the state variable,
and define the spot commodity price process by 
\begin{equation}\label{comm_spot_price}
 S_t:= \exp\left( p^\top Z_t \right) 
\end{equation}
for $ t\in [0,T]$ and $p\in {\mathbb R}^n$. 
This commodity spot price process can be regarded as a flexible extension of the price models proposed in Gibson \& Schwartz (1990) and Schwartz (1997) so to incorporate the skewness of marginal distributions of the log-price process $(\log S_t)_{t\in [0,T]}$.
Indeed, the Gibson \& Schwartz commodity spot price process 
is obtained as follows: 
Let $C=0$, then ${\rm GMSN}(a,\mu_T(Y_0),A,V_T,0)={\mathcal N}(\mu_T(Y_0),V_T)$
and $Z_T=X\sim Y_T$, see Remark 3.1. 
Hence the RMB $(Z_t)$ is equal to the linear Gaussian process
$(Y_t)$. For $n=2$, we consider the parametrisation
\[
d
\begin{pmatrix}
Y^1_t \\ Y^2_t 
\end{pmatrix}
=
\left[
\begin{pmatrix}
\alpha  \\
\kappa \lambda 
\end{pmatrix}
+
\begin{pmatrix}
0 & -1 \\
0 & -\kappa
\end{pmatrix}
\begin{pmatrix}
Y^1_t \\ Y^2_t 
\end{pmatrix}
\right]dt
+
\begin{pmatrix}
\sigma_1& 0 \\
\sigma_2\rho& \sigma_2 \sqrt{1-\rho^2}
\end{pmatrix}
dW_t,
\]
where 
$\alpha,\kappa,\lambda,\sigma_1,\sigma_2>0$ and $-1\le \rho\le 1$.
Then, $S_t = e^{p^\top Y_t}=e^{Y^1_t}$ with $p:=(1,0)^\top$
is precisely the Gibson \& Schwartz (1990) model.
\

By making use of the spot commodity pricing model (\ref{comm_spot_price}), we now have for the forward price process, delivering at time $T$ 
\begin{equation}
 F_{tT}:= {\mathbb E}\left[ S_T\,\vert\, {\mathcal F}^Z_t\right]
={\mathbb E}\left[ e^{p^\top X} \bigm| {\mathcal F}^Z_t\right]
\end{equation}
for $t\in [0,T]$, where the expectation may be taken under the risk-neutral measure, and the filtration $(\F_t^Z)_{t\in[0,T]}$ is generated by the RMB $(Z_t)$.
By Corollary 3.1 and Lemma 3.1, we obtain
\begin{align}\label{pp-int}
F_{tT} =&\int_{{\mathbb R}^n} e^{p^\top x} \pi_t (dx)\nonumber \\
=&\int_{{\mathbb R}^n} e^{p^\top x} f_{\rm GMSN}
\left( x; a,\mu_{T-t}(Z_t),A, V_{T-t},C\right) dx \\
=& F(t,Z_t),\nonumber
\end{align}
where
\begin{multline*}
F(t,z) \\ =\frac{\displaystyle 
\Phi_m\left( C\left( \mu_{T-t}(z)+V_{T-t}p\right)-a;
A+CV_{T-t} C^\top\right)}
{\displaystyle \Phi_m\left( C\mu_{T-t}(z)-a; A+CV_{T-t} C^\top\right)} 
\exp\left[
p^\top \mu_{T-t}(z) +
\frac{1}{2}p^\top V_{T-t}p
\right].
\end{multline*}
On
$(\Omega,{\mathcal F}_{T-}, {\mathbb P}_y, ({\mathcal F}_t^Z)_{t\in [0,T)})$,
we write down the system in stochastic differential form 
\begin{align*}
dF_t =&\nabla_z F(t,Z_t)(\Sigma\Sigma^\top)^{1/2}dU_t, \\
dZ_t =& 
\left[ 
k+ KZ_t
-(\Sigma\Sigma^\top)C^\top
\left( \frac{\nabla \Phi_m}{\Phi_m}\right)
\left( C\mu_{T-t}(Z_t)-a; A+C V_{T-t}C^\top\right)
\right] dt \\
&+(\Sigma\Sigma^\top)^{1/2}dU_t,
\end{align*}
where the SDE determining the dynamics of $(Z_t)$ is given in
Proposition 3.2 (2). 
Furthermore, if we consider the derivative security, whose payoff at the derivative's maturity $T'\le T$ has the form 
$H:=h\left( F_{T'T}\right)$
with some integrable $h:{\mathbb R}_{++} \to {\mathbb R}$, then 
its (arbitrage-free) price at time $t\in [0,T']$ is computed by using 
Proposition 3.2 (1) as follows: 
\begin{eqnarray}
V(t,Z_t)&=&
{\mathbb E}\left[ \e^{-r(T'-t)} h\left( F_{T'T}\right) 
\bigm| {\mathcal F}^Z_t\right]\nn \\
&=&\e^{-r(T'-t)} 
\int_{{\mathbb R}^n} h\left( F( T', y )\right) 
q(t,Z_t; T',y) dy\nn  \\
&=&\e^{-r(T'-t)} 
\int_{{\mathbb R}^n} h\left( F( T', y )\right) 
f_{\rm GMSN}
\left( y; \tilde{a}_t,\mu_{T'-t}(Z_t), \tilde{A}_t,V_{T'-t},
 \tilde{C}_t\right)dy,\nn
\end{eqnarray}
where $r$ is a constant interest rate. In Figure 3, we consider a one-dimensional SNRDB $(Z_t)_{0\le t<T}$
given by (3.1)-(3.3) and (3.10) as the state variable, 
and plot the log-price function 
$L(t,z):=\log F(t,z)$. The parameter values are set to
$n=m=1$, $T=1$, $k=0$, $K=-0.1$, $\Sigma=0.3$, $Z_0=0$, 
$a=1$, $A=1$, $C=-10,-5, 0, 5, 10$, $p=1$, and $t=0.5$. 
We note that $L(t,z)$ is affine 
with respect to $z$ when $C=0$, that is,
\begin{align*}
 L(t,z)=& p\mu_{T-t}(z)+\frac{p^2}{2}V_{T-t} \\
=&
\left\{
\frac{pk}{K}e^{K(T-t)}
+\frac{p^2}{4K}\left( e^{K(T-t)}-1\right)
\right\}
+ p e^{K(T-t)} z.
\end{align*}
The RMB $(Z_t)$ is linear-Gaussian when $C=0$, as we have seen.
We draw the non-linear increasing functions $z \to L(t,z)$ for the values $C\in\{-10, -5,0,5,10\}$. 
\begin{figure}[h]
\begin{center}
\includegraphics[angle=-90, scale=0.5]{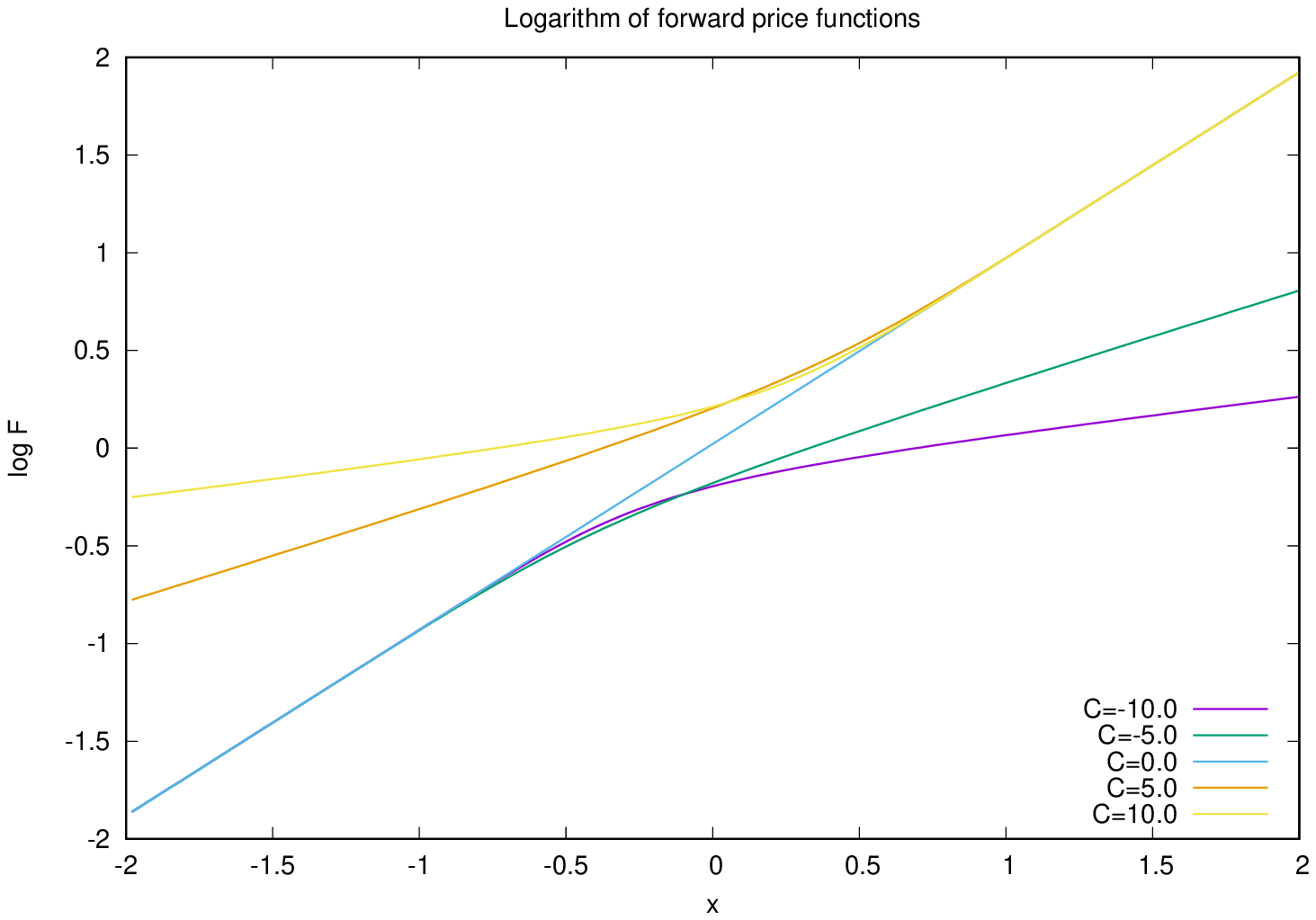}
\caption{Log-forward functions $L(0.5,x)$.}
\end{center}
\end{figure}
Further, in Figure 4, we draw the marginal probability density functions 
\[
 p_{t,T}(x):= \frac{{\mathbb P}\left( F_{tT} \in dx\right)}{dx},
\]
for $t\in\{0.25, 0.5, 0.75, 1\}$, by setting the parameters to
$n=m=1$, $T=1$, $k=0$, $K=-0.1$, $\Sigma=0.3$, $Z_0=0$, 
$a=1$, $A=1$, $C=5$, and $p=1$.
We see that
\[
 p_{t,T}(x)=f_{\rm GMSN}\left( 
F^{-1}(t,x); \tilde{a}_t,\mu_t(Z_0), \tilde{A}_t,V_t,\tilde{C}_t
\right)
\frac{1}{\partial_z F\left( t, F^{-1}(t,x)\right)},
\]
where we use Proposition 3.2 (1), the relation
$F_{tT}=F(t,Z_t)$, and the notation 
$F^{-1}(t,x)$ to denote the inverse function of $F(t,\cdot)$
that satisfies $F\left( t, F^{-1}(t,x)\right)=x$, $x\in {\mathbb R}_{++}$.
\begin{figure}[H]
\begin{center}
\includegraphics[angle=-90, scale=0.5]{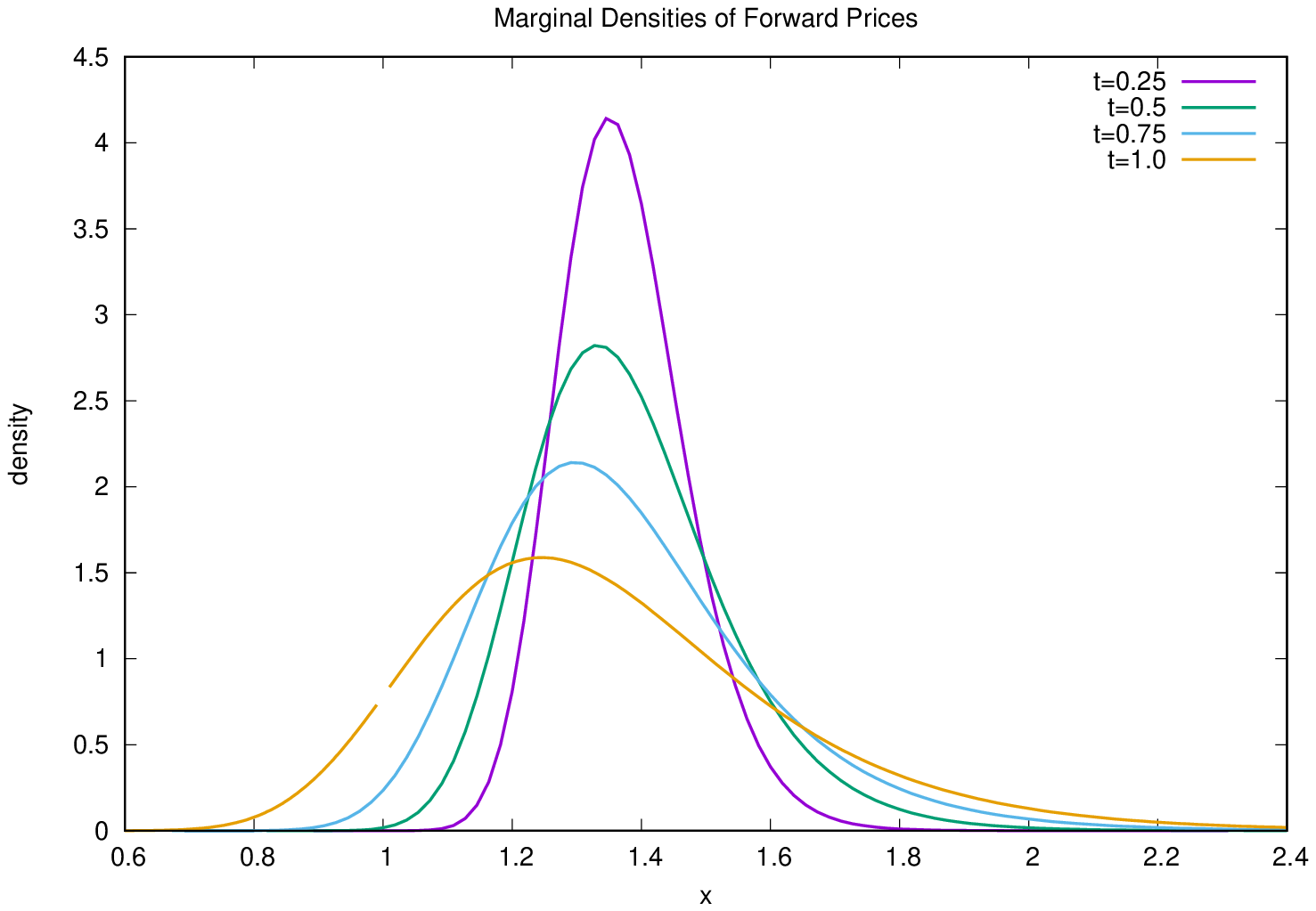}
\caption{Marginal densities of $F_{tT}$.}
\end{center}
\end{figure}
\subsection{Greenhouse gas risk securitisation}
In this section, we associate the RMB to greenhouse gas (GHG) emissions. Here we focus on  man-made emissions (e.g. burning fossil fuel, farming, etc.), which can be affected by varying collective human activity. We regard the man-made emission as a process that in principle can be controlled or on which one can at least intervene. 

We assume that $(Z_t)$ is a randomised Markov bridge with arbitrary marginal distribution at time $t=T$. This RMB is the model for man-made greenhouse gas emissions in the environment. The arbitrary distribution of the `terminal' random variable $Z_T$ is taken to be the target distribution aimed at by, say, an (governmental) emissions policy. Such a target distribution---or indeed the regulatory policy---is specified at $t=0$, when for instance a new (international) law is ratified. Time $T$ is the date by which the target emission level is to be achieved. Because we believe that it is unreasonable to set a fixed target, we opt for a distribution of possible values. We might like to assume that the mean of the distribution coincides with the agreed target value characterising a policy. Such environmental policies come with significant financial commitments. We next present an idea on how environmental policies and related costs may sit at the interface with financial markets and the securitisation of climate risk. 

This is certainly not the first piece of work addressing the question of how financial markets can be involved in reducing GHG emissions. Howison \& Schwarz (2015) provide---in addition to several references---a useful overview in this area of research. One learns that there are mainly three types of reduction mechanisms: (a) the Clean Development Mechanism, (b) the Joint Implementation, and (c) Emissions Trading. The financial mechanism we present next may be seen closest to Emissions Trading in that financial contracts are traded which are written on emission allowance. Carmona \& Hinz (2011) use the so-called risk-neutral valuation approach to price such contracts. While we also use a no-arbitrage approach, though we prefer to work under the real-world measure, we think that the design of the here proposed market instruments written on GHG emissions is different. If we were to look for similarities with existing financial instruments, then the GHG-contracts presented next would be related to inflation-linked or possibly insurance-linked securities, see e.g. Dam et al. (2018) and Barrieu \& Albertini (2009). Such contracts provide financial protection against losses due to, for example, currency depreciation or limit financial liabilities beyond a certain amount. Emission allowance certificates have a different nature: the holder of such a note has paid the right to emit a specified amount of GHG-equivalent. Allowance certificates, as exchanged on e.g. the European Energy Exchange, are not designed to provide financial compensation for potential losses due to (direct or indirect) adverse effects caused by high levels of GHG at a future point in time. While trading emission allowances and the GHG-securities both aim at incentivising a reduction of future GHG levels, the latter mechanism builds in an additional penalisation to the issuer (e.g. a central government)---beyond any specified in international agreements---in a similar way as inflation-linked sovereign bonds result in payments to bond holders by the issuing (governmental) body for failing to keep (future) inflation below a contractually pre-specified level. A governmental GHG policy, as illustrated below, may be viewed as an analogy to a monetary policy and the issuance of GHG-linked securities as an assurance to, e.g., the insurance industry of a central government's plan to reduce GHG levels over a specified time period. We now move one to the stochastic modelling of GHG by use of RMBs and to the no-arbitrage pricing of the proposed GHG-linked securities. 

We consider a central government that has agreed to reducing the GHG amount in the atmosphere by a pre-specified future date. We refer to, e.g., The Paris Agreement (2016)\footnote{The Paris Agreement (2016), United Nations, Framework Convention fro Climate Change, October 2016. http://unfccc.int/paris\_agreement/items/9485.php} whereby the commitment is a reduction of GHG by 80\% by 2050 relative to the GHG amount in 1990, see UK Climate Action Following the Paris Agreement (2016)\footnote{UK Climate Action Following the Paris Agreement (2016), Committee on Climate Change, pp 17-18. https://www.theccc.org.uk/wp-content/uploads/2016/10/UK-climate-action-following-the-Paris-Agreement-Committee-on-Climate-Change-October-2016.pdf}. The target GHG process is modelled by an RMB, where the terminal distribution is a reflection of the uncertainty about the realised target emissions level at time $T$. We suppose that the governing body needs to raise funds to finance the transition to a lower GHG emissions regime by the fixed future date, and aims at doing business with the financial industry. Here we shall consider the insurance industry, whereby (re)insurance firms have a strong interest in a decarbonised future. It is commonly accepted that climate-driven catastrophes are to become more severe (and thus more costly) if GHG emissions do not decrease (c.f. related increase of average global temperature leading to more severe flooding, droughts, hurricanes, etc.). 

We thus consider a governing body that is prepared to sell GHG securities to an insurer where the deal is that (a) if the pre-specified GHG target is met, or even exceeded, then the government retains the funds raised through selling the securities, and (b) if the target is underachieved then the government will have to pay the insurer the difference between the realised GHG amount and the contractually pre-specified target. Such a contract is known as a call option, and in the insurance industry such an instrument may be identified with a stop-loss contract.

We consider the filtered probability space $(\Omega,\F,(\F_t),\PR)$ where $\PR$ is the real probability measure. We take the RMB $(Z_t)$ as a model for the greenhouse gas emissions and introduce a $\PR$-Brownian motion $(B_t)$. We assume $\F_t=\F^B_t\vee\F^Z_t$, where $\F^B_t$ and $\F^Z_t$ are the filtrations generated by $(B_t)$ and $(Z_t)$, respectively. Further, we introduce a stochastic discount factor $(\zeta_t)$ defined by 
\[
\zeta_ t=\exp\left(-\lambda B_t -\tfrac{1}{2}\lambda^2 t - rt\right).
\] 
The interest rate $r$ and the market price of risk $\lambda$ are assumed constant, for simplicity.  Next we introduce a random cash flow $H_T$ at the fixed date $T\ge t$ modelled by $H_T=h(Z_T)$. The no-arbitrage price $(H_t)$, $0\le t\le T$, of the GHG-linked security, with cash flow $H_T$ at time $T$, is given by 
\begin{equation}
H_t = \frac{1}{\zeta_t}\E^{\PR}\left[\zeta_T\, h(Z_T)\,\vert\,\F_t\right].
\end{equation}
Unless one has the view that GHG emissions have a direct impact on the dynamics of the stochastic discount factor $(\zeta_t)$ (and vice versa), we may assume that $\zeta_t$ and $Z_t$ are independent for all $t\in[0,T]$. It then follows that
\[
H_t=\frac{1}{\zeta_t}\E^{\PR}\left[\zeta_T\,\vert\,\F^B_t\right]\E^{\PR}\left[h(Z_T)\,\vert\,\F^Z_t\right],
\]
where it is assumed that 
$h(Z_T)$ is integrable, and further that
\[
H_t=\e^{-r(T-t)}\E\left[h(Z_T)\,\vert\,\F^Z_t\right]
\]
owing to the fact that $\exp(-\lambda B_t -1/2 \lambda^2 t)$ is an
$\{(\F_t),\PR\}$-martingale. The conditional expectation is calculated
in an analogous way to 
Formula (\ref{pp-int}), which gives
\begin{equation}\label{pp-H}
H_t=\e^{-r(T-t)}\int h(x) f_{\rm GMSN}\left(x; a, \mu_{T-t}(Z_t), A, V_{T-t}, C\right) dx,
\end{equation}
Here, in line with the previous section, we assume that $(Z_t)$ is an SNRDB.

{\it Vanilla option}. We next price a stop-loss contract (a call option). We envisage an emissions derivatives market and present price models for GHG emissions options under the said assumptions. We have:
\begin{equation}
V_t	= \e^{-r(T-t)}\E\left[(H_T-K)^+\,\vert\,\F^Z_t\right]=\e^{-r(T-t)}\E\left[\left(h(Z_T)-K\right)^+\,\vert\,\F^Z_t\right]
\end{equation}
where $K$ is a strike price. This equation boils down to calculating
\begin{equation}
V_t=\e^{-r(T-t)}\int \left(h(x) -K\right)^+ 
f_{\rm GMSN}\left(x; a, \mu_{T-t}(Z_t), A, V_{T-t}, C\right) dx.
\end{equation}

This is a first simple example for a GHG-linked option. One might like to consider a richer market where the option maturities are not necessarily associated with an international GHG reduction treaty and dates agreed therein. In fact, a government could issue GHG options with maturities $T_1<T_2<\ldots<T_n$, where $T_n=T$ could coincide with the terminal date (e.g. 2050) by which the GHG reduction must be achieved. The government (and the insurance firms) could then consider holding a basket of GHG options. The policymaker would have goals with shorter terms, which could be more plausible to the buyers in a (financial) insurance market, and the insurers would potentially have access to cash flows after shorter periods of time, which could be used to cover losses. Another type of a contract could be a GHG-linked swap contract, where the insurer exchanges a fixed rate for a floating one that is periodically updated with observations of the realised GHG level. Such a contract gives an opportunity to update the financial liability during a long-dated agreement. 

{\it Swap contract}. Let $T_1<T_2<\ldots<T_i<\ldots<T_n$ be fixed dates. Associated with these dates we have the cash flows $H_{T_i}=g_i(Z_{T_i})$, $i=1,\ldots,n$, where $Z_{T_i}$ has the target distribution of the RMB $(Z^{(i)}_t)_{0\le t\le T_i}$. The integrable functions $g_i$ may differ for different $i$'s. Here the filtration $(\F^Z_t)_{0\le t}$ is generated by the collection of RMBs $\{(Z^{(i)}_t)_{0\le t\le T_i}\}_{i=1,\ldots,n}$. Then we consider the following swap-like structure with price process $(Sw_t)$:
\begin{equation}
Sw_t=\E\left[\sum^n_i\e^{-r(T_i-t)}\left(H_{T_i}-K\right)\,\big\vert\,\F^Z_t\right]
\end{equation}
where $K$ is a strike value and $(\F^Z_t)$ is generated by the RMBs $\{(Z^{(1)}_t), (Z^{(2)}_t),\ldots, (Z^{(n)}_t)\}$. The expression for the price $Sw_t$ at time $t$ can be written as
\begin{equation}
Sw_t=\sum^n_i\e^{-r(T_i-t)}\E\left[g_i(Z_{T_i})\,\vert\,\F^Z_t\right]-K\sum^n_i \e^{-r(T_i-t)}.
\end{equation}
By setting $Sw_t=0$ we obtain the so-called {\it fair swap value} $K_t$ at time $t$ given by
\begin{eqnarray}
K_t&=&\sum^n_i\E\left[g_i(Z_{T_i})\,\vert\,\F_t\right]\big/\sum^n_i\e^{-r(T_i-t)}\nn\\
     &=&\sum^n_i\int g_i(x_i)
f_{\rm GMSN}\left(x_i; a_i, \mu^{(i)}_{T_i-t}(Z^{(i)}_t), A_i, V^{(i)}_{T_i-t}, C_i\right) dx_i \big/\sum^n_i\e^{-r(T_i-t)},\nn\\
\end{eqnarray}
where here $\{(Z^{(i)}_t)\}_{i=1,\ldots,n}$ are assumed to be SNRDBs. It is $K_t$ that seller (government) and buyer (insurer) trade at. In the case of a swap, the value of $K_t$ is directly related to how likely a government will actually meet the GHG target (or even reduce to below the target). Also, since $Sw_t$ may take positive and negative values, the trade counterparts in principle rebalance their accounts as they go along in contrast to call option contracts where there would be no cashflows for a potentially long period of time.

\begin{figure}[H]
\centering
\includegraphics[scale=.55]{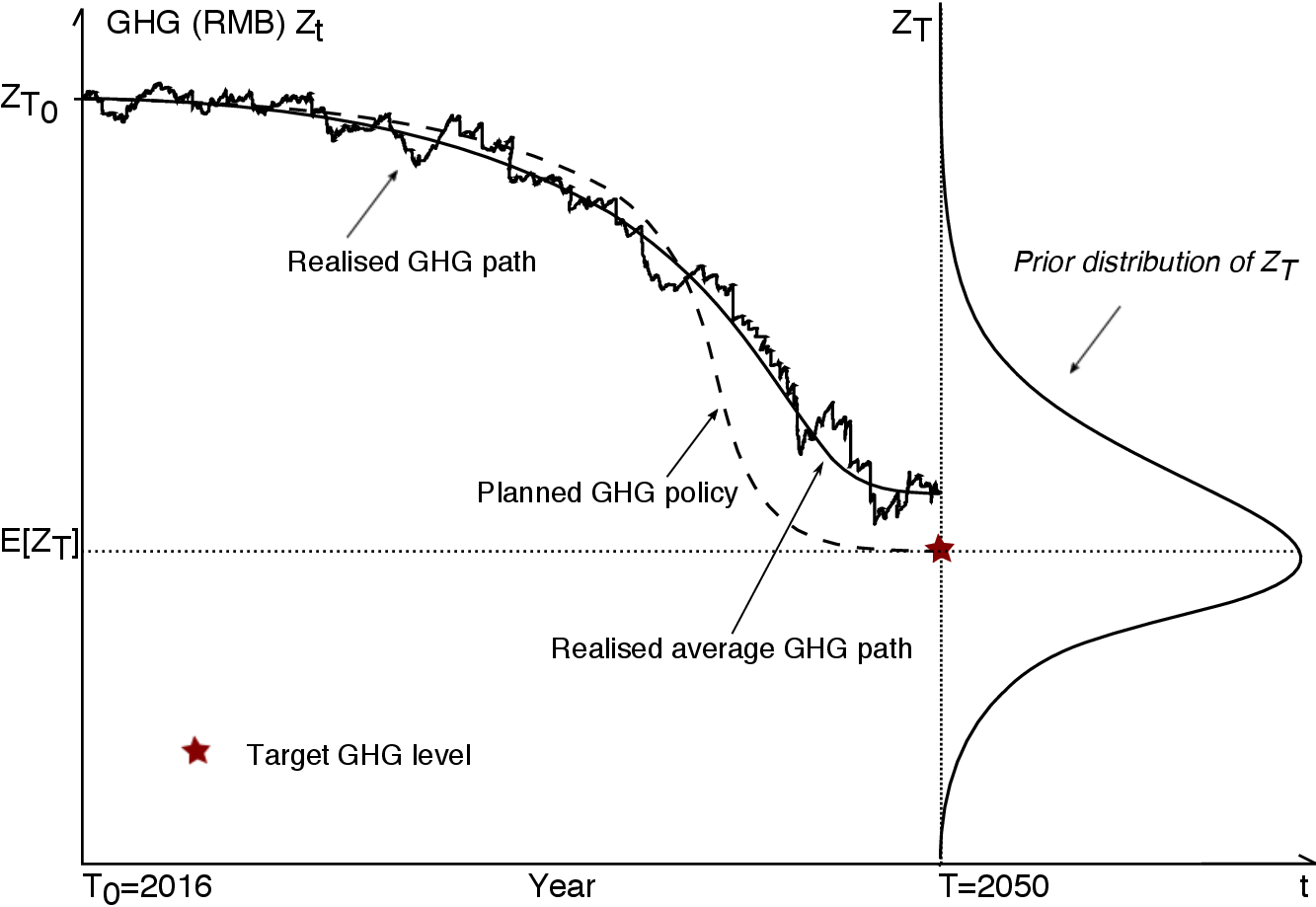}
\caption{Greenhouse gas (GHG) level modelled by an RMB $(Z_t)$.}\label{GHG-RMB}
\end{figure}

In Figure \ref{GHG-RMB}, we describe how an RMB may be applied to model GHG emission. At $T_0$, $Z_{T_0}$ denotes the observed initial level of emitted GHG. Also at $T_0$, the prior distribution of the random GHG level $Z_T$ is specified. The mean $\E[Z_T]$ of the prior distribution is set to match the target GHG level that the guided emissions policy is supposed to meet by the pre-specified date $T$. The dashed line depicts the GHG emission plan a policymaker commits to and wishes to implement. The stochastic path is one possible realisation of the GHG process $(Z_t)$, which is observed as one moves forward in time . The solid line is the implicit drift of the process $(Z_t)$, which in the depicted example does not match the planned GHG policy. At the end of the policy period, that is at $T$, the terminal GHG level $Z_T$ exceeds the pre-specified target value, which the policymaker had committed to realise. We emphasize that at first the realised GHG level (and in particular the realised drift) is below the dashed GHG policy line, thus indicating, in this particular example, that the implemented emission policy achieved a stronger abatement than intend up to that point in time.
We conclude the application of RMBs to GHG-risk securitization by taking a closer look at what model parameters a policymaker may influence. In other words, what parameters in the chosen RMB may capture trends a policymaker may wish to control during the planned GHG abatement period. A candidate parameter is $k$ in the linear Gaussian diffusion (\ref{Ydiff}), underlying the dynamics of $(Z_t)$, because it affects the drift of the RMB. This can be seen explicitly via Corollary \ref{Cor} and the set of equations (\ref{meanvar}). 

The idea is that the GHG level increases or declines depending on what level of $k$ is selected. The parameter $k$, in turn, affects the price dynamics of a financial security written on $(Z_t)$. Other degrees of freedom, such as the {\it a priori} target distribution of $Z_T$, also have an impact, but here we shall focus on the sensitivity of $(H_t)$ on $k$. As seen in Formula (\ref{pp-H}), the value $(H_t)$ at time $t$ is impacted by $k$ via $\mu_{T-t}$ in the GMSN-density. We express the `$k$-sensitivity' of $H_t$ at time $t\in[0,T)$ by $\partial_k H_t$. In the following proposition, we also give an explicit expression for the so-called `Delta-sensitivity' $\partial_z H_t$. This is the measure of how much the price $H_t$ at time $t$ varies given a change in the underlying value of $Z_t$ at time $t$. This result is given for the case that the SNRDB is multivariate process.
\newpage
\begin{prop}
Consider the SNRDB, given by (3.1)--(3.3) and (3.10).
Let $h:=(h_1,\dots,h_d)^\top:{\mathbb R}^n \to {\mathbb R}^d$, 
assuming that $h(X)$ is integrable. 
It then holds that
\begin{eqnarray}
\partial_{z} {\mathbb E}\left[ h(X) \bigm| {\mathcal F}^Z_t\right]
&=& {\rm Cov}\left[ h(X), X \bigm| {\mathcal F}^Z_t\right]^\top
 V^{-1}_{T-t}e^{(T-t)K}, \\
\partial_k {\mathbb E}\left[ h(X) \bigm| {\mathcal F}^Z_t\right]
&=&
{\rm Cov}\left[ h(X), X \bigm| {\mathcal F}^Z_t\right]^\top
V^{-1}_{T-t}\int_0^{T-t} e^{sK}ds,
\end{eqnarray}
where we use (2.17), (2.8), $\pi_t(h)=\Pi(t,Z_t;h)$, 
and write
\begin{align*}
\partial_{z} {\mathbb E}\left[ h(X) \bigm| {\mathcal
 F}^Z_t\right]&:=\partial_z {\Pi}(t,Z_t; h)
=\left( \partial_{z_j} {\Pi}(t,Z_t;h_i)\right)_{1\le i\le d, 1\le j\le n},
\\
\partial_{k} {\mathbb E}\left[ h(X) \bigm| {\mathcal
 F}^Z_t\right]&:=
\left( \partial_{k_j}{\Pi}(t,Z_t;h_i)\right)_{1\le i\le d, 1\le j\le n}.
\end{align*}
\end{prop} 
\begin{proof}
Equation (4.13) is obtained from Proposition 2.3 and 
(2.15)--(2.17) in Remark 2.8.
To show (4.14), we note that, in the 
right-hand side of the expression, 
\begin{equation}
{\mathbb E}\left[ h(X) \bigm| {\mathcal F}^Z_t\right]
= \int_{{\mathbb R}^n} h(x) 
f_{\rm GMSN}\left(x; a,\mu_{T-t}(Z_t),A,V_{T-t},C \right)dx,
\end{equation}
where we use Corollary 3.1, 
only $\mu_{T-t}(Z_t)$ depends on $Z_t$ and $k$.
We see that
\begin{eqnarray}
 \partial_z \mu_{T-t}(Z_t)&=&e^{(T-t)K}, \\
 \partial_k \mu_{T-t}(Z_t)&=&\int_0^{T-t} e^{sK}ds. 
\end{eqnarray} 
From (4.13), (4.15), (4.16) and the chain rule of differentiation,
we deduce
\begin{equation}
\partial_b 
\left\{
\int_{{\mathbb R}^n} h(x) 
f_{\rm GMSN}\left(x; a,b,A,V_{T-t},C \right)dx\right\} 
\Bigm|_{b=\mu_{T-t}(Z_t)}
={\rm Cov}\left[ h(X),X | {\mathcal F}^Z_t\right] V_{T-t}^{-1}.
\end{equation}
From (4.15), (4.17) and (4.18), 
we see that
\begin{align*}
\partial_k {\mathbb E}\left[ h(X) \bigm| {\mathcal F}^Z_t\right]
=& 
\partial_b
\left\{
\int_{{\mathbb R}^n} h(x) 
f_{\rm GMSN}\left(x; a,b,A,V_{T-t},C \right)dx
\right\} \Bigm|_{b=\mu_{T-t}(Z_t)} 
 \partial_k \mu_{T-t}(Z_t) \\
=&{\rm Cov}\left[ h(X),X | {\mathcal F}^Z_t\right] V_{T-t}^{-1}
\int_0^{T-t} e^{sK}ds.
\end{align*}
\end{proof}
\begin{rem}
Letting $h(z)=z$ in Proposition 4.1, 
we see that 
\begin{eqnarray*}
\partial_z {\mathbb E}\left[ X \bigm| {\mathcal F}^Z_t\right]
&=&
{\rm Var}\left[ X \bigm| {\mathcal F}^Z_t\right]
V^{-1}_{T-t}e^{(T-t)K}\ge 0, \\
\partial_k {\mathbb E}\left[ X \bigm| {\mathcal F}^Z_t\right]
&=&
{\rm Var}\left[ X \bigm| {\mathcal F}^Z_t\right]
V^{-1}_{T-t}\int_0^{T-t} e^{sK}ds\ge 0.
\end{eqnarray*}
\end{rem}
\section*{Acknowledgements}
A large part of this article was written while J. Sekine was visiting 
the Department of Mathematics, University College London. 
UCL is acknowledged for their kind hospitality. A. Macrina thanks Osaka University and The Daiwa Anglo-Japanese Foundation for financial support provided for this research collaboration. The authors are grateful to S. N. Cohen, C. A. Garcia Trillos, M. R. Grasselli and G. W. Peters for interesting discussions and suggestions, and to anonymous reviewers for feedback leading to improvements of this paper. A first version of this work, entitled Filtering with Randomised Markov Bridges, appeared on \url{https://arxiv.org/abs/1411.1214v1}.


\begin{thebibliography}{99}

\bibitem{AC1}
{\sc Azzalini, A. and A. Capitanio} (1999)
Statistical applications of the multivariate skew normal distribution.
{\it J. Roy. Statist. Soc.} B {\bf 61}, 579--602.

\bibitem{AC2}
{\sc Azzalini, A. and A. Capitanio} (2013)
``The Skew-Normal and Related Families.
Institute of Mathematical Statistics Monographs.'' 
{\it Cambridge University Press}. 

\bibitem{AD}
{\sc Azzalini, A. and A. Dalla Valle} (1996)
The multivariate skew-normal distribution.
{\it Biometrika},  {\bf 83}(4), 715--726.

\bibitem{BA}
{\sc Barrieu, P. and L. Albertini} (2009)
The Handbook of Insurance-Linked Securities.
{\it Wiley Finance}. 

 \bibitem{B}
{\sc Baudoin, F.} (2002).
Conditioned stochastic differential equations:
theory, examples and application to finance.
{\it Stochastic Processes and their Applications}, 
{\bf 100}, 109--145. 

 \bibitem{BHM}
{\sc Brody, D. C., L. P. Hughston and A. Macrina} (2008).
Information-based asset pricing.
{\it International Journal of Theoretical and Applied Finance},
{\bf 11}(1), 107--142.

%

 \bibitem{CD}
{\sc Carmona, R. and Hinz, J.} (2011).
Risk-Neutral Models for Emission Allowance Prices and Option Valuation. 
{\it Management Science}, {\bf 57}(8),
1453--1488.

 \bibitem{CD}
{\sc \c{C}etin, U. and Danilova, A.} (2016).
Markov bridges: SDE representation. 
{\it Stochastic Processes and their Applications}, {\bf 126}(3),
651--679.

 \bibitem{CUB}
{\sc Chaumont, L, and Uribe Bravo, G.} (2011).
Markovian bridges: weak continuity and pathwise constructions.
{\it The Annals of Probability},
{\bf 39}(2), 609--647.


 \bibitem{CMNS}
{\sc Cr\'epey, S., A. Macrina, M. T. Nguyen and D. Skovmand} (2016).
Rational Multi-Curve Models with Counterparty-Risk Adjustments.
{\it Quantitative Finance} {\bf 16}(6), 847-866.

\bibitem{Dam}
{\sc Dam, H., A. Macrina, D. Sloth and D. Skovmand} (2018). Rational model for inflation-linked derivatives. \url{http://dx.doi.org/10.2139/ssrn.3110835} 

 \bibitem{FHM}
{\sc Filipovi\'c, D., L. P. Hughston and A. Macrina} (2012).
Conditional density models for asset pricing.
{\it International Journal of Theoretical and Applied Finance},
{\bf 15}(1),  1250002, 24 pp.
 
 \bibitem{FPM}
{\sc Fitzsimmons, P., J. Pitman, and M. Yor} (1993).
Markovian bridges: construction, palm interpretation, and splicing. 
{\it Seminar on Stochastic Processes, 1992.
(Progress in Probability, Volume 33, 
Editors: E. \c{C}inlar, K. L. Chung, M. J. Sharpe, R. F. Bass, and K. Burdzy),}
101--134, Birkh\"auser.

 \bibitem{GS}
{\sc Gibson, R., and E. S. Schwartz} (1990). 
Stochastic convenience yield and the pricing of oil contingent claims.
{\it The Journal of Finance} {\bf 45}, 959--976. 

 \bibitem{GGD}
{\sc Gupta, A. K., G. Gonz\'alez-Far\'ias, and J. A.  Dom\'inguez-Molina}
(2004).
A multivariate skew normal distribution.
{\it Journal of Multivariate Analysis}, 89, 181-190.

 \bibitem{HHM}
{\sc Hoyle, E., L. P. Hughston and A. Macrina} (2011).
L\'evy random bridges and the modelling of
financial information.
{\it Stochastic Processes and their Applications}, 
{\bf 121},
856--884.

 \bibitem{HHM15}
{\sc Hoyle, E., L. P. Hughston and A. Macrina} (2015).
Stable-1/2 Bridges and Insurance.
{\it Advances in Mathematics of Finance}, Banach Center Publications, Volume 104, Institute of Mathematics, Polish Academy of Sciences, Warsaw. 95--120

 \bibitem{HS}
 {\sc Howison, S., D. Schwarz} (2015). Risk-Neutral Pricing of Financial Instruments in Emission Markets: A Structural Approach. {\it SIAM Review}, {\bf 57}(1), 95--127.


 \bibitem{Macrina}
 {\sc Macrina, A.} (2014). Heat Kernel Models for Asset Pricing. {\it International Journal of Theoretical and Applied Finance}, {\bf 17}(7), 1450048, 34 pp.

 \bibitem{P}
{\sc Protter, P.} (2005). 
Stochastic Integration and Differential Equations (Second Edition). 
{\it Springer}.

 \bibitem{S}
{\sc Schwartz, E. S.} (1997).
The Stochastic Behavior of Commodity Prices: 
Implications for Valuation and Hedging.
{\it The Journal of Finance}, {\bf 52}(3), 923--973.

\end{thebibliography}
\end{document}